\documentclass[11pt,a4paper]{amsart}
\usepackage{amsmath,amsthm,amssymb}
\usepackage[all]{xy}

\DeclareMathOperator{\pol}{pol}

\newcommand{\typI}{\boldsymbol1}%{\mathbf{I}}
\newcommand{\typII}{\boldsymbol2}%{\mathbf{II}}
\newcommand{\typIII}{\boldsymbol3}%{\mathbf{III}}
\newcommand{\typIV}{\boldsymbol4}%{\mathbf{IV}}
\newcommand{\typV}{\boldsymbol5}%{\mathbf{V}}
\newcommand{\typVI}{\boldsymbol6}%{\mathbf{VI}}
\newcommand{\typVII}{\boldsymbol7}%{\mathbf{VII}}
\newcommand{\typVIII}{\boldsymbol8}%{\mathbf{VIII}}
\newcommand{\typIX}{\boldsymbol9}%{\mathbf{IX}}
\swapnumbers
\theoremstyle{definition}
\newtheorem{point}{}[section]

\theoremstyle{plain}

\newtheorem{theorem}[point]{Theorem}
\newtheorem{state}[point]{Statement}

\newcommand{\marginextend}[1]{ \addtolength{\oddsidemargin}{-#1}  \addtolength{\evensidemargin}{-#1}\addtolength{\textwidth}{#1}\addtolength{\textwidth}{#1}}
\newcommand{\updownextend}[1]{ \addtolength{\topmargin}{-#1}  \addtolength{\textheight}{#1}
\addtolength{\textheight}{#1}}
\DeclareMathOperator{\Id}{Id}

\marginextend{1.25cm}
\updownextend{0.00cm}\allowdisplaybreaks[4]
\title[FQ operations II. Examples, monoaxiality, scalings]{Fermionic quantum operations: a computational framework II.
Examples, monoaxiality, scalings}
\author{Gyula Lakos}
\email{lakos@cs.elte.hu}
\address{Department of Geometry, E\"otv\"os Lor\'and University, P\'azm\'any P\'eter s.~1/C,  Budapest, H--1117, Hungary}
\begin{document}
\maketitle
\begin{abstract}
The objective of this series of papers is to recover information regarding the behaviour of FQ operations
in the case $n=2$, and  FQ conform-operations in the case $n=3$.
In this second part we show some arithmetically constructible examples of FQ operations ($n=2$),
concentrating on  monoaxiality, related extensions,  and (hyper)scaling.
\end{abstract}
\section*{Introduction}
Advancing analytical tools for noncommuting sets of operators has a long history.
In contrast to the relatively simple case of single matrices cf. \cite{DS}, \cite{Rin},  where the spectrum plays a decisive role,
the study of noncommuting sets of operators is much more difficult.
Systematic approaches bring forth various organization principles, like
the use of the joint spectrum \cite{Tay}, powerful applications of perturbation expansions \cite{KVV} (also see references
therein), operator orderings \cite{NSS} (with much prehistory in quantum physics), application of Clifford algebras  \cite{BDS}
along with  efforts to integrate various approaches, cf.~\cite{Jef}.
The study of concrete operator functions is also developing, especially operator means, see e.~g.  \cite{HP},
and in preserver problems, which is a sort of an invariant theory for noncommutative systems, cf. \cite{P}.
Combining Clifford systems and free analysis,
we intend to further the idea of using Clifford systems as base systems, whose perturbations can be investigated:
More specifically,
FQ operations were introduced in \cite{L1} as some kind of quasi linear algebraic tools for non-commuting operators.
These can be approached on global (analytic) and local / punctual (formal) level.
The objective of this series of papers (cf. \cite{L2}) is to recover information regarding the behaviour of FQ operations
in the case $n=2$, and  FQ conform-operations in the case $n=3$.
In this second part we show some arithmetically constructible examples of FQ operations ($n=2$),
concentrating on  monoaxiality, related extensions,  and (hyper)scaling.
Here, in contrast to \cite{L2}, we consider only natural (conjugation-invariant) FQ operations with sign-linearity.
Moreover, we aim here to orthogonal invariant FQ operations.
The outline of the paper is as follows.
Section \ref{sec1} contains some arithmetical constructions of FQ operations.
In particular, we construct an FQ orthogonalization procedure.
Due to this arithmeticity (as opposed to constructions by iterative methods)
examples here will be rather simple; and  some of them are
 arguably the simplest ones.
The emphasis is on showing up some ``major'' types.
There are many similar ones with various tradeoffs in their behaviour.
Some of them  have less than perfect properties (like no involutivity or idempotence) but
they might serve as initial objects in other constructions.
In Section \ref{sec2} we consider the axial extension procedure.
In Section \ref{sec3} we consider  hyperscaling conditions and
see why they are too strong.

\section{Some arithmetically constructible examples}\label{sec1}
In this section we will consider certain  analytic FQ operations.
Our computations will be valid when the involved terms ($\mathcal A_C$, $\mathcal B^{-1/2}$, etc.) are well-defined,
but always valid on the formal domain.
In this section all FQ operations will be symmetric, and, in fact, orthogonal invariant.
Whenever we consider an operation we will give its first-order expansion terms relative to  the mixed base of \cite{L2},
which tells much about the character of the operation.
Beyond that, in order to distinguish the various cases, we will not give  terms of higher order,
which are hard to interpret, but we  indicate their scalar scaling properties (i. e. possible scaling variables in the mixed base).
Fortunately, scalar scaling properties are rather easy to check or disprove in concrete cases.
\begin{point}
In if $A,B\in\mathfrak A$, and the straight segment $\{(1-t)A+tB\,:\,t\in[0,1]\}$ lies in $\mathfrak A^\star$, then one defines the
geometric mean
\[(A\star B)^{1/2}:=\int_{t=0}^{2\pi} \frac{1}{A^{-1}\cos^2 t+B^{-1}\sin^2 t}\,\frac{\mathrm dt}{2\pi},\]
which is also the common value
 \[A(A^{-1}B)^{1/2}=B(B^{-1}A)^{1/2}=(AB^{-1})^{1/2}B=(BA^{-1})^{1/2}A.\]

One can see that if $A$ and $B$ are skew-involutions, then $(A\star B)^{1/2}$ is also a skew-involution, furthermore,
\[(A\star B)^{1/2}=\pol\tfrac12(A+B);\]
where we have used the notation $\pol X=X(-X^2)^{-1/2}$.
Moreover, in that case, $A$ and $B$ are conjugates of each other by $(A\star B)^{1/2}$.
We will occasionally use the notation $|A|=(-A^2)^{1/2}$.
\end{point}

\begin{point}
We define the pseudoscalar FQ operations
left axis $\mathcal A_L$ by
\[\mathcal A_L(A_1,A_2):=\pol -A_1A_2^{-1}=\pol A_2A_1^{-1};\]
right axis  $\mathcal A_R$
\[\mathcal A_R(A_1,A_2):=\pol A_2^{-1}A_1=\pol -A_1^{-1}A_2;\]
and (central) axis $\mathcal A_C$ by
 \[\mathcal A_C(A_1,A_2):=(\mathcal A_L(A_1,A_2)\star\mathcal A_R(A_1,A_2))^{1/2}.\]
Furthermore we define the scalar FQ operation biaxiality  $\mathcal B$ as
\[\mathcal B(A_1,A_2):= (-\mathcal A_L(A_1,A_2)\cdot\mathcal A_R(A_1,A_2))^{1/2} ;\]
in this case
\[\mathcal B(A_1,A_2)^{-1}= (-\mathcal A_R(A_1,A_2)\cdot\mathcal A_L(A_1,A_2))^{1/2} .\]

One can see that the corresponding terms in their (formal) expansion are
\[\mathcal A_L:\qquad \hat {\mathbf P}^{[12]}_0=[1]\qquad
\hat {\mathbf P}^{[12]}_1=\begin{bmatrix}0&2&0&0&0&0&2&0\end{bmatrix};\]
\[\mathcal A_R:\qquad \hat {\mathbf P}^{[12]}_0=[1]\qquad
\hat {\mathbf P}^{[12]}_1=\begin{bmatrix}0&-2&0&0&0&0&2&0\end{bmatrix};\]
\[\mathcal A_C:\qquad \hat {\mathbf P}^{[12]}_0=[1]\qquad
\hat {\mathbf P}^{[12]}_1=\begin{bmatrix}0&0&0&0&0&0&2&0\end{bmatrix};\]
and
\[\mathcal B:\qquad \hat {\mathbf P}^{[0]}_0=[1]\qquad
\hat {\mathbf P}^{[0]}_1=\left[ \begin {array}{cccccccc} 0&2&0&0&0&0&0&0\end {array}\right];\]
\[\mathcal B^{-1}:\qquad \hat {\mathbf P}^{[0]}_0=[1]\qquad
\hat {\mathbf P}^{[0]}_1=\left[ \begin {array}{cccccccc} 0&-2&0&0&0&0&0&0\end {array}\right].\]
\end{point}

It is immediate but useful to keep in mind that, using the notation $A_1=\Id_1$, $A_2=\Id_2$,
\begin{enumerate}
\item$\mathcal A_L\mathcal A_1=- A_1\mathcal A_R,\quad A_L\mathcal A_2=- A_2\mathcal A_R$;
\item
$\mathcal A_L=\mathcal B^{\pm1}\mathcal A_L\mathcal B^{\pm1}=\mathcal B\mathcal A_C=\mathcal A_C\mathcal B^{-1}=\mathcal B^2\mathcal A_R=
\mathcal B\mathcal A_R\mathcal B^{-1}=\mathcal A_R\mathcal B^{-2}=-\mathcal A_{C}\mathcal A_{R}\mathcal A_{C}$;
\item
$\mathcal A_R=\mathcal B^{\pm1}\mathcal A_R\mathcal B^{\pm1}=\mathcal B^{-1}\mathcal A_C=\mathcal A_C\mathcal B=\mathcal B^{-2}\mathcal A_L=
\mathcal B^{-1}\mathcal A_L\mathcal B=\mathcal A_L\mathcal B^2=-\mathcal A_{C}\mathcal A_{L}\mathcal A_{C}$;
\item $\mathcal A_C=\mathcal B^{\pm1}\mathcal A_C\mathcal B^{\pm1}=\mathcal B\mathcal A_R=\mathcal A_R\mathcal B^{-1}=
\mathcal B^{-1}\mathcal A_L=\mathcal A_L\mathcal B$;
\item $\mathcal B=-\mathcal A_{L}\mathcal B^{-1}\mathcal A_{L}=
-\mathcal A_{C}\mathcal B^{-1}\mathcal A_{C}=-\mathcal A_{R}\mathcal B^{-1}\mathcal A_{R}
=-\mathcal A_{C}\mathcal A_{R}=-\mathcal A_{L}\mathcal A_{C}$;
\item $\mathcal B^{-1}=-\mathcal A_{L}\mathcal B\mathcal A_{L}=-\mathcal A_{C}\mathcal B\mathcal A_{C}=-\mathcal A_{R}\mathcal B\mathcal A_{R}
=-\mathcal A_{R}\mathcal A_{C}=-\mathcal A_{C}\mathcal A_{L}$.
\end{enumerate}
One can also see that $\mathcal A_L,\mathcal A_C,\mathcal A_R,\mathcal B,\mathcal B^{-1}$
satisfy scaling invariances in  $\hat r_1,\hat r_2,\hat r_3,\hat r_4,\hat r_5$.
\begin{point}Taking multiplication positionwise,
$\mathcal A_L(A_1,A_2)\cdot (A_1,A_2)=-(A_1,A_2)\cdot \mathcal A_R(A_1,A_2)$.
This allows us to define the axial turn operation $\mathcal T$ as
\[\mathcal T(A_1,A_2):= \mathcal A_L(A_1,A_2)\cdot (-A_2, A_1)=(A_2, -A_1)\cdot\mathcal A_R(A_1,A_2).\]
The corresponding terms in its expansion are
\[\mathcal T:\qquad \hat {\mathbf P}^{[1]}_0=[1]\qquad
\hat {\mathbf P}^{[1]}_1=\left[ \begin {array}{cccccccc} 1&1&1&-1&-1&1&1&1\end {array}\right].\]

One can see
\[\mathcal O^{\mathrm{fSy}}(A_1,A_2)=\frac12((A_1,A_2)+\mathcal T(A_1,A_2)),\]
where $\mathcal O^{\mathrm{fSy}}$ is the canonical realization of the symmetric conform-orthogonalization procedure from \cite{L1}.
The corresponding terms in its expansion are
\[\mathcal O^{\mathrm{fSy}}:\qquad \hat {\mathbf P}^{[1]}_0=[1]\qquad
\hat {\mathbf P}^{[1]}_1=\left[ \begin {array}{cccccccc} 1&1&1&0&0&1&1&1\end {array}\right].\]
We define the anticonform orthogonalization operation ${\mathcal O}^{\mathrm{afSy}}$ by
\[{\mathcal O}^{\mathrm{afSy}}(A_1,A_2):= (-\mathcal O^{\mathrm{fSy}}(A_1,A_2)_1^{-1},-\mathcal O^{\mathrm{fSy}}(A_1,A_2)_2^{-1}).\]
${\mathcal O}^{\mathrm{afSy}}$ also produces a floating Clifford system, except this operation is not bivariant but antivariant.
(${\mathcal O}^{\mathrm{afSy}}$ is set up so that it is Clifford conservative.)
${\mathcal O}^{\mathrm{afSy}}$ can also be realized by a closed integral formula, cf.  \cite{L1}.
The corresponding terms in its expansion are
\[{\mathcal O}^{\mathrm{afSy}}:\qquad \hat {\mathbf P}^{[1]}_0=[1]\qquad
\hat {\mathbf P}^{[1]}_1=\left[ \begin {array}{cccccccc} -1&-1&-1&0&0&1&1&1\end {array}\right].\]

One can check that
the following hold:
\begin{enumerate}
\item
$\mathcal A_L \circ {\mathcal T}=\mathcal A_L$,\quad
$\mathcal A_C \circ {\mathcal T}=\mathcal A_C$,\quad
$\mathcal A_R \circ {\mathcal T}=\mathcal A_R$,\quad
$\mathcal B \circ {\mathcal T}=\mathcal B$;
\item
$\mathcal A_L \circ {\mathcal O}^{\mathrm{fSy}}=\mathcal A_L$,\quad
$\mathcal A_C \circ {\mathcal O}^{\mathrm{fSy}}=\mathcal A_C$,\quad
$\mathcal A_R \circ {\mathcal O}^{\mathrm{fSy}}=\mathcal A_R$,\quad
$\mathcal B \circ {\mathcal O}^{\mathrm{fSy}}=\mathcal B$;
\item
$\mathcal A_L \circ {\mathcal O}^{\mathrm{afSy}}=\mathcal A_R$,\quad
$\mathcal A_C \circ {\mathcal O}^{\mathrm{afSy}}=\mathcal A_C$,\quad
$\mathcal A_R \circ {\mathcal O}^{\mathrm{afSy}}=\mathcal A_L$,\quad
$\mathcal B \circ {\mathcal O}^{\mathrm{afSy}}=\mathcal B^{-1}$.
\end{enumerate}

One can also see that $\mathcal T$ satisfies scaling invariances in variables $\hat r_1,\hat r_2,\hat r_3,\hat r_4,\hat r_5$;
${\mathcal O}^{\mathrm{fSy}}$ satisfies scaling invariances in variables $\hat r_1,\hat r_2,\hat r_3,\tilde r_5$;
${\mathcal O}^{\mathrm{afSy}}$ satisfies scaling invariances in variables $\hat r_1,\hat r_2,\hat r_3,\tilde r_4$.
(Here $\tilde r_4$ and $\tilde r_5$ refer to some mixed scaling conditions in $\hat r_4 / \hat r_5$.)
\end{point}

\begin{point}
We say that the pair $(A_1,A_2)$ is monoaxial if $\mathcal A_L(A_1,A_2)$ and $\mathcal A_R(A_1,A_2)$,
and hence, $\mathcal A_C(A_1,A_2)$, are equal to each other.
This happens if and only if $\mathcal B(A_1,A_2)=1$.
Another equivalent formulation is that $\mathcal A_C(A_1,A_2)$ anticommutes with $A_1,A_2$.
Indeed, the identites $A_1\mathcal A_R(A_1,A_2)A_1^{-1}=\mathcal A_L(A_1,A_2)$ and $A_2\mathcal A_R(A_1,A_2)A_2^{-1}=\mathcal A_L(A_1,A_2)$
imply that in the monoaxial case $A_1$ and $A_2$ commute with the axes, i. e. with the axis.
Conversely, if $\mathcal A_C(A_1,A_2)$ commutes with $A_1$, $A_2$ then it commutes with $\mathcal A_L(A_1,A_2)=\pol -A_1A_2^{-1}$;
on the other hand, as a general rule, we know that  $\mathcal A_L(A_1,A_2)$ conjugated  by $\mathcal A_C(A_1,A_2)$  is $\mathcal A_R(A_1,A_2)$;
consequently $\mathcal A_L(A_1,A_2)=\mathcal A_R(A_1,A_2)$.
\end{point}
\begin{theorem}
If $(A_1,A_2)$ is monoaxial with axis $\mathcal A_C=\mathcal A_C(A_1,A_2)$, and
\[1+\tilde r_3^{\mathrm{x}}:=\left(-\frac{A_1^2+A_2^2-A_1\mathcal A_CA_2+A_2\mathcal A_CA_1}{4}\right)^{1/2}\]
exists;  then, we claim,
\[Q_1=\pol\,\, \mathcal O^{\mathrm{fSy}}(A_1,A_2)_1,\qquad Q_2=\pol\,\, \mathcal O^{\mathrm{fSy}}(A_1,A_2)_2\]
forms a Clifford system; and with $\tilde r_4^{\mathrm{x}}:=-A_1Q_1+Q_2A_2$,  the deomposition
\[(A_1,A_2)=((1+\tilde r_3^{\mathrm{x}}+\tilde r_4^{\mathrm{x}})Q_1,(1+\tilde r_3^{\mathrm{x}}-\tilde r_4^{\mathrm{x}})Q_2)\]
is valid. (Actually, this is the circular decomposition with respect to $(Q_1,Q_2)$.)
Here $Q_1,Q_2$ commute with $ \tilde r_3^{\mathrm{x}}$ and anticommute with $\mathcal A_C$ and $\tilde r_4^{\mathrm{x}}$; moreover, $\mathcal A_C=Q_1Q_2$.

Furthermore,
\[\mathcal O^{\mathrm{fSy}}(A_1,A_2)=((1+\tilde r_3^{\mathrm{x}})Q_1,(1+\tilde r_3^{\mathrm{x}})Q_2),\]
and
\[\mathcal O_{\mathrm{m}}^{\mathrm{Sy}}(A_1,A_2)=(Q_1,Q_2);\]
where the latter equation is understood such that it provides an analytical realization of
$\underline{\mathcal O}^{\mathrm{Sy}}$ (for some monoaxial elements, though).
\begin{proof}
From monoaxiality one can deduce
\[(\mathcal O^{\mathrm{fSy}}(A_1,A_2)_1)^2=(\mathcal O^{\mathrm{fSy}}(A_1,A_2)_2)^2=
\tfrac14\left({A_1^2+A_2^2-A_1\mathcal A_CA_2+A_2\mathcal A_CA_1}\right). \]
After that, it is a straightforward computation.
\end{proof}
\end{theorem}
\begin{point}
We say that a vectorial FQ operation is monoaxial if its result is always monoaxial (e.~g.: FQ orthogonalizations).
We define the vectorial FQ operations left monoaxialization
\[\mathcal M_L(A_1,A_2):= (A_1,A_2)\cdot\mathcal B(A_1,A_2)^{-1} ;\]
right monoaxialization,
\[\mathcal M_R(A_1,A_2):=\mathcal B(A_1,A_2)^{-1} \cdot(A_1,A_2);\]
and central monoaxialization
\[\mathcal M_C(A_1,A_2):=\mathcal B(A_1,A_2)^{-1/2} \cdot(A_1,A_2)\cdot \mathcal B(A_1,A_2)^{-1/2}.\]
These are conjugates of each other. It is easy to see that
$\mathcal M_L$, $\mathcal M_R$, $\mathcal M_C$ are monoaxial
with axes $\mathcal A_L$, $\mathcal A_R$, $\mathcal A_C$, respectively;
and they act trivially on monoaxial pairs.

The corresponding terms in their (formal) expansion are
\[\mathcal M_L:\qquad \hat {\mathbf P}^{[1]}_0=[1]\qquad
\hat {\mathbf P}^{[1]}_1=\left[ \begin {array}{cccccccc} -1&1&1&1&1&1&1&1\end {array} \right] ;\]
\[\mathcal M_R:\qquad \hat {\mathbf P}^{[1]}_0=[1]\qquad
\hat {\mathbf P}^{[1]}_1=\left[ \begin {array}{cccccccc} 1&-1&1&1&1&1&1&1\end {array} \right] ;\]
\[\mathcal M_C:\qquad \hat {\mathbf P}^{[1]}_0=[1]\qquad
\hat {\mathbf P}^{[1]}_1=\left[ \begin {array}{cccccccc} 0&0&1&1&1&1&1&1\end {array} \right] .\]

We  define the axial conjugation  operation $\mathcal C$ as
\[\mathcal C(A_1,A_2):= \mathcal A_C(A_1,A_2)\cdot (A_1, A_2)\cdot\mathcal A_C(A_1,A_2)=
\mathcal B(A_1,A_2)^{-1}\cdot (A_1, A_2)\cdot\mathcal B(A_1,A_2)^{-1}.\]
In terms of its expansion
\[\mathcal C:\qquad \hat {\mathbf P}^{[1]}_0=[1]\qquad
\hat {\mathbf P}^{[1]}_1=\left[ \begin {array}{cccccccc} -1&-1&1&1&1&1&1&1\end {array}\right].\]

One can easily check that with the choices $X,Y=L,C,R$ and conventions $-L=R, -C=C, -R=L$,
the following hold:
\begin{enumerate}
\item$\mathcal A_Y \circ \mathcal M_X=\mathcal A_{X},\quad\mathcal A_{-X} \circ \mathcal C=\mathcal A_{X} $;
\item$\mathcal B \circ \mathcal M_X=1,\quad \mathcal B \circ \mathcal C=\mathcal B^{-1}$;
\item${\mathcal T} \circ \mathcal M_X=\mathcal M_{X}\circ {\mathcal T},\quad
{\mathcal T} \circ \mathcal C=\mathcal C\circ {\mathcal T}$;
\item${\mathcal O}^{\mathrm{fSy}} \circ \mathcal M_X=\mathcal M_{X}\circ {\mathcal O}^{\mathrm{fSy}},\quad
{\mathcal O}^{\mathrm{fSy}} \circ \mathcal C=\mathcal C\circ {\mathcal O}^{\mathrm{fSy}}$;
\item${\mathcal O}^{\mathrm{afSy}} \circ \mathcal M_X=\mathcal M_{-X}\circ {\mathcal O}^{\mathrm{afSy}},\quad
{\mathcal O}^{\mathrm{afSy}} \circ \mathcal C=\mathcal C\circ {\mathcal O}^{\mathrm{afSy}}$;
\item$\mathcal M_Y \circ \mathcal M_X=\mathcal M_{-X} \circ \mathcal C=\mathcal C \circ \mathcal M_X=\mathcal M_X,
\quad\mathcal  C \circ \mathcal C=\Id$;
\item $\mathcal B(A_1,A_2)=1 \Leftrightarrow  \mathcal M_X(A_1,A_2)=(A_1,A_2).$
\end{enumerate}
Furthermore, $\mathcal M_L,\mathcal M_C,\mathcal M_R,\mathcal C$
satisfy scaling invariances in  $\hat r_1,\hat r_2,\hat r_3,\hat r_4,\hat r_5$.
\end{point}
\begin{point}
At this point one can easily define an orthogonalization procedure by
\[{\mathcal O}^{\mathrm{mSy}}:={\mathcal O}^{\mathrm{Sy}}_{\mathrm{m}}\circ\mathcal M_C.\]
Indeed, after monoaxialization, ${\mathcal O}^{\mathrm{Sy}}_{\mathrm{m}}$ can be applied.
In its expansion, the corresponding terms of first-order are
\[{\mathcal O}^{\mathrm{mSy}}:\qquad \hat {\mathbf P}^{[1]}_0=[1]\qquad
\hat {\mathbf P}^{[1]}_1=\left[ \begin {array}{cccccccc} 0&0&0&0&0&1&1&1\end {array}\right].\]
From the definition it is easy to see that
\begin{enumerate}
\item$\mathcal A_Y \circ {\mathcal O}^{\mathrm{mSy}}=\mathcal A_{C},\quad\mathcal B \circ {\mathcal O}^{\mathrm{mSy}}=1 $;
\item${\mathcal O}^{\mathrm{mSy}}
={\mathcal O}^{\mathrm{mSy}}\circ{\mathcal O}^{\mathrm{fSy}}
={\mathcal O}^{\mathrm{mSy}}\circ{\mathcal O}^{\mathrm{afSy}}
={\mathcal O}^{\mathrm{mSy}}\circ\mathcal T
={\mathcal O}^{\mathrm{mSy}}\circ\mathcal M_C
={\mathcal O}^{\mathrm{mSy}}\circ\mathcal C.$
\end{enumerate}

Now, ${\mathcal O}^{\mathrm{mSy}}$ has the same expansion in first order as
 $\underline{\mathcal O}^{\mathrm{Sy}}$, but they already differ in  second order.
More informatively, one can see that (relative to the mixed base) $\underline{\mathcal O}^{\mathrm{Sy}}$ has a single scaling invariance,
in the variable $\hat r_3$, i. e. scalar homogeneity.
On the other hand, one can show that ${\mathcal O}^{\mathrm{mSy}}$ has  scaling invariances
in the variables $\hat r_1,\hat r_2,\hat r_3$ (but it does not satisfy metric trace commutativity then, of course).
\end{point}
\begin{point}
In order to have a more explicit form, we define the
 scalar FQ operation axial length (squared)  $\mathcal L^{\mathrm{fx}}$ by
\[\mathcal L^{\mathrm{fx}}(A_1,A_2):=\frac14      \left(-A_1\mathcal B^{-1}A_1-A_2\mathcal B^{-1}A_2+A_1\mathcal A_CA_2-A_2\mathcal A_CA_1\right), \]
axial volume  $\mathcal V^{\mathrm{fm}}$ by
\[\mathcal V^{\mathrm{fm}}(A_1,A_2):=\mathcal B^{-1/2} \mathcal L^{\mathrm{fx}}(A_1,A_2)\mathcal B^{-1/2}, \]
and the axial pseudodeterminant $\mathcal D^{\mathrm{fx}}$ as
\[\mathcal D^{\mathrm{fx}}(A_1,A_2):=\frac14 \left(A_1\mathcal A_CA_1+A_2\mathcal A_CA_2+A_1\mathcal B^{-1}A_2-A_2\mathcal B^{-1}A_1\right) .    \]
where $\mathcal B=\mathcal B(A_1,A_2)$, and $\mathcal A_C=\mathcal A_C(A_1,A_2)$.
Terms in their expansion are
\[\mathcal L^{\mathrm{fx}}:\qquad \hat {\mathbf P}^{[0]}_0=[1]\qquad
\hat {\mathbf P}^{[0]}_1=\left[ \begin {array}{cccccccc} 0&2&2&0&0&0&0&0\end {array} \right] ;\]
\[\mathcal V^{\mathrm{fm}}:\qquad \hat {\mathbf P}^{[0]}_0=[1]\qquad
\hat {\mathbf P}^{[0]}_1=\left[ \begin {array}{cccccccc} 0&0&2&0&0&0&0&0\end {array} \right] ;\]
\[\mathcal D^{\mathrm{fx}}:\qquad \hat {\mathbf P}^{[12]}_0=[1]\qquad
\hat {\mathbf P}^{[12]}_1=\left[ \begin {array}{cccccccc} 0&0&2&0&0&0&2&0\end {array} \right] .\]
One can easily check that
the following hold:
\begin{enumerate}
\item $\mathcal D^{\mathrm{fx}}=\mathcal A_L\mathcal L^{\mathrm{fx}}=\mathcal L^{\mathrm{fx}}\mathcal A_R$,\,
$\mathcal L^{\mathrm{fx}}\mathcal B^{-1}= -\mathcal D^{\mathrm{fx}} \mathcal A_C,\quad
\mathcal B^{-1}\mathcal L^{\mathrm{fx}}= -\mathcal A_C\mathcal D^{\mathrm{fx}} $;

\item
$\mathcal L^{\mathrm{fx}}\circ {\mathcal O}^{\mathrm{fSy}}=\mathcal L^{\mathrm{fx}}$,
$\mathcal L^{\mathrm{fx}}\circ {\mathcal O}^{\mathrm{mSy}}=1$,
$\mathcal L^{\mathrm{fx}}\circ {\mathcal O}^{\mathrm{afSy}}=(\mathcal L^{\mathrm{fx}})^{-1}$;

\item
$\mathcal V^{\mathrm{fm}}\circ {\mathcal O}^{\mathrm{fSy}}=\mathcal V^{\mathrm{fm}}$,
$\mathcal V^{\mathrm{fm}}\circ {\mathcal O}^{\mathrm{mSy}}=1$,
$\mathcal V^{\mathrm{fm}}\circ {\mathcal O}^{\mathrm{afSy}}=(\mathcal V^{\mathrm{fm}})^{-1}$;

\item
$\mathcal D^{\mathrm{fx}}\circ {\mathcal O}^{\mathrm{fSy}}=\mathcal D^{\mathrm{fx}}$,
$\mathcal D^{\mathrm{fx}}\circ {\mathcal O}^{\mathrm{mSy}}=\mathcal A_C$,
$\mathcal D^{\mathrm{fx}}\circ {\mathcal O}^{\mathrm{afSy}}=-(\mathcal D^{\mathrm{fx}})^{-1}$;

\item${\mathcal O}^{\mathrm{mSy}}=(\mathcal V^{\mathrm{fm}})^{-1/2}\mathcal B^{-1/2}{\mathcal O}^{\mathrm{fSy}}\mathcal B^{-1/2}
=\mathcal B^{-1/2}{\mathcal O}^{\mathrm{fSy}}\mathcal B^{-1/2}(\mathcal V^{\mathrm{fm}})^{-1/2}$

$=\mathcal B^{-1/2}(\mathcal L^{\mathrm{fx}}\mathcal B^{-1})^{-1/2}{\mathcal O}^{\mathrm{fSy}}\mathcal B^{-1/2}
=\mathcal B^{-1/2}{\mathcal O}^{\mathrm{fSy}}(\mathcal B^{-1}\mathcal L^{\mathrm{fx}})^{-1/2}\mathcal B^{-1/2}$.
\end{enumerate}
One can also prove that $\mathcal L^{\mathrm{fx}},\mathcal V^{\mathrm{fm}},\mathcal D^{\mathrm{fx}}$ are scaling invariant in variables $\hat r_1,\hat r_2,\hat r_3$.

Notice that $\mathcal L^{\mathrm{fx}}$ and $\mathcal D^{\mathrm{fx}}$ has some variants given by
\[\mathcal L^{\mathrm{x} }(A_1,A_2):=\frac12      \left(-A_1\mathcal B^{-1}A_1-A_2\mathcal B^{-1}A_2\right), \]
\[\mathcal L^{\mathrm{cx}}(A_1,A_2):=\frac12      \left(A_1\mathcal A_CA_2-A_2\mathcal A_CA_1\right), \]
\[\mathcal D^{\mathrm{x} }(A_1,A_2):=\frac12 \left(A_1\mathcal B^{-1}A_2-A_2\mathcal B^{-1}A_1\right)=
\mathcal A_L\cdot\mathcal L^{\mathrm{cx} }(A_1,A_2)=\mathcal L^{\mathrm{cx} }(A_1,A_2)\cdot\mathcal A_R,    \]
\[\mathcal D^{\mathrm{cx} }(A_1,A_2):=\frac12 \left(A_1\mathcal A_CA_1+A_2\mathcal A_CA_2\right)=
\mathcal A_L\cdot\mathcal L^{\mathrm{x} }(A_1,A_2)=\mathcal L^{\mathrm{x} }(A_1,A_2)\cdot\mathcal A_R ,    \]
where  $\mathcal A_L=\mathcal A_L(A_1,A_2)$, $\mathcal A_R=\mathcal A_R(A_1,A_2)$;
which can be related to the even simpler operations
\[\mathcal L(A_1,A_2):=\frac12      \left(-A_1^2-A_2^2\right), \]
\[\mathcal D(A_1,A_2):=\frac12 \left(A_1A_2-A_2A_1\right) .    \]

Here $\mathcal L^{\mathrm{x} }$ and $\mathcal D^{\mathrm{cx} }$  satisfy scalings in variables $\hat r_1,\hat r_2,\hat r_3,\hat r_5$;
$\mathcal L^{\mathrm{cx} }$ and $\mathcal D^{\mathrm{x} }$ satisfy scalings in variables $\hat r_1,\hat r_2,\hat r_3,\hat r_4$;
$\mathcal L$ satisfies scalings in variables $\hat r_2,\hat r_3,\hat r_5$ and an $\hat r_1/\hat r_4$ mixed condition;
$\mathcal D$ satisfies scalings in variables $\hat r_2,\hat r_3,\hat r_4$ and an $\hat r_1/\hat r_5$ mixed condition
(not detailed here).

One can consider the polarizations $\pol\mathcal D^{\mathrm{fx}},\pol\mathcal D^{\mathrm{x}},\pol\mathcal D^{\mathrm{cx}}$ and $\mathcal A_D:=\pol \mathcal D$.
They have the same first-order expansions as $\mathcal A_C$.
$\pol\mathcal D^{\mathrm{fx}}$ scales in $\hat r_1,\hat r_2,\hat r_3$; $\pol\mathcal D^{\mathrm{x}}$ scales in $\hat r_1,\hat r_2,\hat r_3,\hat r_4$;
$\pol\mathcal D^{\mathrm{cx}}$ scales in $\hat r_1,\hat r_2,\hat r_3,\hat r_5$;
$\mathcal A_D=\pol\mathcal D$ scales in $\hat r_2,\hat r_3,\hat r_4$.
Or, one can take the ``volumes'' $|\mathcal D^{\mathrm{fx}}|,|\mathcal D^{\mathrm{x}}|,|\mathcal D^{\mathrm{cx}}|$ and $\mathcal V:=|\mathcal D|$.
They have the same first-order expansions as $\mathcal V^{\mathrm{fm}}$.
$|\mathcal D^{\mathrm{fx}}|$ scales in $\hat r_1,\hat r_2,\hat r_3$; $\pol\mathcal D^{\mathrm{x}}$ scales in $\hat r_1,\hat r_2,\hat r_3,\hat r_4$;
$|\mathcal D^{\mathrm{cx}}|$ scales in $\hat r_1,\hat r_2,\hat r_3,\hat r_5$;
$\mathcal V=|\mathcal D|$ scales in $\hat r_2,\hat r_3,\hat r_4$.
(Remark: $\mathcal V^{\mathrm{fm}}\neq|\mathcal D^{\mathrm{fx}}|$, they differ in third order, or, more qualitatively,
$\mathcal V^{\mathrm{fm}}$ commutes with $\mathcal A_C$ multiplicatively, but $|\mathcal D^{\mathrm{fx}}|$ does not.)
\end{point}
\begin{point}
We can define the vectorial FQ operation axial unitalization $\mathcal U^{\mathrm{fx}}$ by
\[\mathcal U^{\mathrm{fx}}(A_1,A_2):= \mathcal  B^{1/2}(\mathcal V^{\mathrm{fm}})^{-1/4} \mathcal  B^{-1/2} \cdot (A_1,A_2) \cdot
\mathcal   B^{-1/2}(\mathcal V^{\mathrm{fm}})^{-1/4}\mathcal  B^{1/2} ; \]
and the vectorial FQ operation axial amplitude inversion $\mathcal K^{\mathrm{fx}}$ by
\[\mathcal K^{\mathrm{fx}}(A_1,A_2):= \mathcal   B^{1/2} (\mathcal V^{\mathrm{fm}})^{-1/2}\mathcal  B^{-1/2} \cdot (A_1,A_2)
\cdot \mathcal  B^{-1/2}(\mathcal V^{\mathrm{fm}})^{-1/2}\mathcal   B^{1/2} ; \]
where $\mathcal B=\mathcal B(A_1,A_2)$, $\mathcal V^{\mathrm{fm}}=\mathcal V^{\mathrm{fm}}(A_1,A_2)$.
Terms in their expansion are
\[\mathcal U^{\mathrm{fx}}:\qquad \hat {\mathbf P}^{[1]}_0=[1]\qquad
\hat {\mathbf P}^{[1]}_1=\left[ \begin {array}{cccccccc} 1&1&0&1&1&1&1&1\end {array} \right] ;\]
\[\mathcal K^{\mathrm{fx}}:\qquad \hat {\mathbf P}^{[1]}_0=[1]\qquad
\hat {\mathbf P}^{[1]}_1=\left[ \begin {array}{cccccccc} 1&1&-1&1&1&1&1&1\end {array} \right] .\]

One can see that
\begin{enumerate}\item
$\mathcal K^{\mathrm{fx}} \circ \mathcal K^{\mathrm{fx}}=\Id$,
$\mathcal U^{\mathrm{fx}} \circ \mathcal U^{\mathrm{fx}}=
\mathcal U^{\mathrm{fx}} \circ \mathcal K^{\mathrm{fx}}=\mathcal K^{\mathrm{fx}} \circ \mathcal U^{\mathrm{fx}}
=\mathcal U^{\mathrm{fx}}$;

\item
$\mathcal A_L\circ\mathcal U^{\mathrm{fx}}=\mathcal A_L$,
$\mathcal A_C\circ\mathcal U^{\mathrm{fx}}=\mathcal A_C$,
$\mathcal A_R\circ\mathcal U^{\mathrm{fx}}=\mathcal A_R$,
$\mathcal B\circ\mathcal U^{\mathrm{fx}}=\mathcal B$;

\item
$\mathcal A_L\circ\mathcal K^{\mathrm{fx}}=\mathcal A_L$,
$\mathcal A_C\circ\mathcal K^{\mathrm{fx}}=\mathcal A_C$,
$\mathcal A_R\circ\mathcal K^{\mathrm{fx}}=\mathcal A_R$,
$\mathcal B\circ\mathcal K^{\mathrm{fx}}=\mathcal B$;

\item $\mathcal U^{\mathrm{fx}}$ and $\mathcal K^{\mathrm{fx}}$ commutes with $\mathcal M_L,\mathcal M_C,\mathcal M_R,\mathcal C,
{\mathcal O}^{\mathrm{fSy}},{\mathcal O}^{\mathrm{afSy}},\mathcal T$ with respect to $\circ$;

\item $ {\mathcal O}^{\mathrm{mSy}}=\mathcal U^{\mathrm{fx}}\circ{\mathcal O}^{\mathrm{fSy}}\circ\mathcal M_C
=\mathcal U^{\mathrm{fx}}\circ{\mathcal O}^{\mathrm{afSy}}\circ\mathcal M_C$ ,
$\mathcal K^{\mathrm{fx}}\circ{\mathcal O}^{\mathrm{fSy}}\circ\mathcal M_C={\mathcal O}^{\mathrm{afSy}}\circ\mathcal M_C$;

\item $\mathcal L^{\mathrm{fx}}\circ\mathcal U^{\mathrm{fx}}=\mathcal B$, $\mathcal D^{\mathrm{fx}}\circ\mathcal U^{\mathrm{fx}}=\mathcal A_C$;

\item
$\mathcal V^{\mathrm{fm}} \circ \mathcal K^{\mathrm{fx}}=(\mathcal V^{\mathrm{fm}})^{-1}$, $ \mathcal V^{\mathrm{fm}} \circ \mathcal U^{\mathrm{fx}}=1$; $
\mathcal V^{\mathrm{fm}}(A_1,A_2)=1 \Leftrightarrow  \mathcal U^{\mathrm{fx}}(A_1,A_2)=(A_1,A_2)$.
\end{enumerate}

Furthermore, $\mathcal U^{\mathrm{fx}}$ and $\mathcal K^{\mathrm{fx}}$  has scaling invariances in variables $\hat r_1,\hat r_2,\hat r_3$.
\end{point}
\begin{point}
We see that we have three commuting sets of actions $\mathcal M_C,\mathcal C$  and $\mathcal U^{\mathrm{fx}},\mathcal K^{\mathrm{fx}}$ and $\mathcal O^{\mathrm{fSy}},\mathcal T$;
acting principally in variables $\hat r_1,\hat r_2$ and $\hat r_3$ and $\hat r_4,\hat r_5$ respectively.
We can take the composition $\mathcal O^{\mathrm{mSy}}=\mathcal M_C\circ\mathcal U^{\mathrm{fx}}\circ\mathcal O^{\mathrm{fSy}}$ which we have already seen.
But one can also define
\[\mathcal I_C^{\mathrm{fy}}:=\mathcal C\circ\mathcal K^{\mathrm{fx}}\circ\mathcal T.\]
It has expansion terms
\[\mathcal I_C^{\mathrm{fy}}:\qquad \hat {\mathbf P}^{[1]}_0=[1]\qquad
\hat {\mathbf P}^{[1]}_1=\left[ \begin {array}{cccccccc} -1&-1&-1&-1&-1&1&1&1\end {array} \right] .\]
It is an involution,
\[(\mathcal I_C^{\mathrm{fy}})^2=\Id.\]
Here $\mathcal I_C^{\mathrm{fy}}$ has scaling invariances in variables $\hat r_1,\hat r_2,\hat r_3$.
One can check that
\[\mathcal I_C^{\mathrm{fy}}(A_1,A_2)=(\mathcal B(A_1,A_2)\star \mathcal L^{\mathrm{fx}}(A_1,A_2))^{-1/2}\cdot
\mathcal T(A_1,A_2) \cdot(\mathcal B(A_1,A_2)\star \mathcal L^{\mathrm{fx}}(A_1,A_2))^{-1/2}.\]
One can also define
\[\mathcal I_L^{\mathrm{fy}}(A_1,A_2):=\mathcal B(A_1,A_2)^{-1}\cdot \mathcal T(A_1,A_2)\cdot \mathcal L^{\mathrm{fx}}(A_1,A_2)^{-1},\]
\[\mathcal I_R^{\mathrm{fy}}(A_1,A_2):=\mathcal L^{\mathrm{fx}}(A_1,A_2)^{-1}\cdot \mathcal T(A_1,A_2)\cdot\mathcal B(A_1,A_2)^{-1}.\]
These operations have the  same expansion in first order, and the same scaling properties, and they are also involutive.
It is easy to check that
\begin{enumerate}
\item
$\mathcal L^{\mathrm{fx}} \circ\mathcal I_C^{\mathrm{fy}}=\mathcal L^{\mathrm{fx}} \circ\mathcal I_L^{\mathrm{fy}}=
\mathcal L^{\mathrm{fx}} \circ\mathcal I_R^{\mathrm{fy}}=(\mathcal L^{\mathrm{fx}} )^{-1}$;
\item
$\mathcal D^{\mathrm{fx}} \circ\mathcal I_C^{\mathrm{fy}}=\mathcal D^{\mathrm{fx}} \circ\mathcal I_L^{\mathrm{fy}}=
\mathcal D^{\mathrm{fx}} \circ\mathcal I_R^{\mathrm{fy}}=-(\mathcal D^{\mathrm{fx}} )^{-1}$.
\end{enumerate}
There are some variants given by
\[\mathcal I_L^{\mathrm{y}}(A_1,A_2):=\mathcal B(A_1,A_2)^{-1}\cdot \mathcal T(A_1,A_2)\cdot \mathcal L^{\mathrm{cx}}(A_1,A_2)^{-1},\]
\[\mathcal I_R^{\mathrm{y}}(A_1,A_2):=\mathcal L^{\mathrm{cx}}(A_1,A_2)^{-1}\cdot \mathcal T(A_1,A_2)\cdot\mathcal B(A_1,A_2)^{-1},\]
\[\mathcal I_L^{\mathrm{cy}}(A_1,A_2):=\mathcal B(A_1,A_2)^{-1}\cdot \mathcal T(A_1,A_2)\cdot \mathcal L^{\mathrm{x}}(A_1,A_2)^{-1},\]
\[\mathcal I_R^{\mathrm{cy}}(A_1,A_2):=\mathcal L^{\mathrm{x}}(A_1,A_2)^{-1}\cdot \mathcal T(A_1,A_2)\cdot\mathcal B(A_1,A_2)^{-1}.\]
They have the same first-order expansions as above;
$\mathcal I_L^{\mathrm{y}},\mathcal I_R^{\mathrm{y}}$
have scalings in  $\hat r_1,\hat r_2,\hat r_3,\hat r_4$, and
$\mathcal I_L^{\mathrm{cy}},\mathcal I_R^{\mathrm{cy}}$
have scalings in  $\hat r_1,\hat r_2,\hat r_3,\hat r_5$; but these variants are not involutive. Also,
\begin{enumerate}
\item[(3)] $\mathcal I_C^{\mathrm{fy}},\mathcal I_L^{\mathrm{fy}},\mathcal I_R^{\mathrm{fy}}
,\mathcal I_L^{\mathrm{y}},\mathcal I_R^{\mathrm{y}},\mathcal I_L^{\mathrm{cy}},\mathcal I_R^{\mathrm{cy}}$
commute with $\mathcal T$ and $\mathcal C$ with respect to composition.
\end{enumerate}\end{point}

$\mathcal I_C^{\mathrm{y}}$ is reasonably nice but it lacks $\hat r_4$ scaling and hence affine invariance.
What would be more interesting is to have a similar involutive and transposition invariant operation but with affine invariance properties
(i.~e.~with orthogonal invariance and $\hat r_3$ and $\hat r_4$ scalings).
\section{The axial extension procedure}\label{sec2}
\begin{point}\label{po:mono}
Some of the latter examples are examples ``axial extensions'' which we explain as follows.
Suppose that $\Xi$ is a (sort of) conjugation-invariant FQ operation  defined  at least on monoaxial elements.
Then we can set
\[\Xi^{\mathrm{x}}(A_1,A_2):=\Xi(A_1\mathcal B^{-1},A_2\mathcal B^{-1})\cdot\mathcal B
=\mathcal B\cdot\Xi(\mathcal B^{-1}A_1,\mathcal B^{-1}A_2)\]\[=\mathcal B^{1/2}\cdot\Xi(\mathcal B^{-1/2}A_1\mathcal B^{-1/2},\mathcal B^{-1/2}A_2\mathcal B^{-1/2})\cdot\mathcal B^{1/2}\]
where $\mathcal B=\mathcal B(A_1,A_2)$. The equality of the various presentations is a consequence of conjugation invariance.
For and FQ operation, being axial (extension) can be understood as having bivariance with respect to the biaxility operation.
Now, it is easy to see that this extension (or modification) process acts conservatively on monoaxial elements.

This axial extension procedure has a variant defined by
\[\not\!\Xi^{\mathrm{cx}}(A_1,A_2):=\pm\Xi(A_1\mathcal A_C,A_2\mathcal A_C)\cdot\mathcal A_C^{-1}
=\pm\mathcal A_C^{-1}\cdot\Xi(\mathcal A_C\mathcal A_1,\mathcal A_CA_2),\]
where  $\mathcal A_C=\mathcal A_C(A_1,A_2)$, with $\pm=+$ in the scalar and vectorial cases, and $\pm=-$ in the pseudscalar case.
Equivalence follows from naturality.
This is more or the less the same as the original version:
Formally, or analytically, if the environment is appropriate the following hold:
If $\Xi$ is scalar, then $\not\!\Xi^{\mathrm{cx}}=\mathcal A_L\cdot \Xi^{\mathrm{x}}=\Xi^{\mathrm{x}}\cdot \mathcal A_R$;
if $\Xi$ is vectorial, then $\not\!\Xi^{\mathrm{cx}}=\Xi^{\mathrm{x}}$;
if $\Xi$ is pseudoscalar, then $\not\!\Xi^{\mathrm{cx}}=-\mathcal A_L\cdot \Xi^{\mathrm{x}}=-\Xi^{\mathrm{x}}\cdot \mathcal A_R$.
(It interchanges scalar and pseudoscalar.)
\end{point}
\begin{point}
The axial extension procedure can be understood on the formal level as follows.
What happens is that for perturbations $(A_1,A_2)=(Q_1+R_1,Q_2+R_2)$ of $(Q_1,Q_2)$ we have decomposition
\begin{multline*}(Q_1+R_1,Q_2+R_2)=((1+\hat r_1+\hat r_2+\hat r_3+\hat r_4+\hat r_5+\hat r_6+\hat r_7+\hat r_8)Q_1,\\
(1+\hat r_1-\hat r_2+\hat r_3-\hat r_4-\hat r_5+\hat r_6+\hat r_7-\hat r_8)Q_2).\end{multline*}
Now, if $(A_1,A_2)$, and its axis is $Q_1Q_2$, then
according to what we have seen, monoaxility makes $(A_1,A_2)$ anticommute with $Q_1Q_2$; hence
$\hat r_1=\hat r_2=\hat r_7=\hat r_8=0$ in the decomposition. So,
\[(Q_1+R_1,Q_2+R_2)=((1+\hat r_3+\hat r_4+\hat r_5+\hat r_6)Q_1,
(1+\hat r_3-\hat r_4-\hat r_5+\hat r_6)Q_2).\]
Ultimately, one has a restricted formal FQ calculus for monoaxial operations, which is similar
to the original one but with a restricted set of variables  $\hat r_3,\hat r_4,\hat r_5,\hat r_6$.
In this case $\mathcal O^{\mathrm{fSy}}(A_1,A_2)=(Q_1,Q_2)$ if and only if $\hat r_6=0$.
Hence, there is a similar discussion of conjugation-invariance leading to the eliminability of $6$-indices
(as opposed to  $ 6, 7, 8 $-indices). One can also see that it can also be described by
the (restricted version) of the $6$-hyperscaling condition.

Notice that using the usual expansion formulas we cannot create non-monoaxial pairs from monoaxial pairs because an expression of
$\hat r_3,\hat r_4,\hat r_5,\hat r_6$ still commutes with $Q_1Q_2$. (``No escape from monoaxiality.'')
Similarly, if one has a pseudoscalar Clifford conservative operation, then it must be multiplicatively
conjugate to $\mathcal A_C$, hence by the reason above, it must be $\mathcal A_C$ itself.
Also notice that in the monoaxial calculus, multiplication by $\mathcal A_C$ makes a bijective correspondence between scalar and
pseudoscalar operations.

\begin{theorem}\label{theo:ax}
In terms of expansions regarding the mixed base, being an axial extension means
that lower $\{1,2,7,8\}$-indices can be eliminated from $\hat p^{[s]}_{\iota_1,\ldots,\iota_r}$;
or, taking conjugation-invariance into account, it means the eliminability of $\{1,2\}$-indices, hence reduction to lower
$\{3,4,5\}$-indices.
\begin{proof}
This follows from conservativity on monoaxial elements and using any of the extension formulas from \ref{po:mono}.
\end{proof}
\end{theorem}
In the Clifford conservative case, for an axial extension $\Xi^{\mathrm{x}}$,
in the scalar case $\hat q^{[0]}_1=0$ and $\hat q^{[0]}_2=2$;
in the vectorial case $\hat q^{[1]}_1=\hat q^{[1]}_2=1$;
in the pseudoscalar case $\hat q^{[12]}_1=\hat q^{[12]}_2=0$.

 We can also interpret $\Xi^{\mathrm m}:=\Xi\circ\mathcal M_C$ as a kind of axial extension with
$\hat q^{[s]}_1=\hat q^{[s]}_2=0$. Or, we can consider
\[\Xi^{\mathrm{y}}(A_1,A_2):=\mathcal B^{-1}\cdot\Xi(A_1\mathcal B^{-1},A_2\mathcal B^{-1})
=\Xi(\mathcal B^{-1}A_1,\mathcal B^{-1}A_2)\cdot\mathcal B^{-1}\]
\[=\mathcal B^{-1/2}\cdot\Xi(\mathcal B^{-1/2}A_1\mathcal B^{-1/2},\mathcal B^{-1/2}A_2\mathcal B^{-1/2})\cdot\mathcal B^{-1/2}\]
as an axial extension. Then, for the extension $\Xi^{\mathrm y}$ of a Clifford conservative operation $\Xi$,
in the scalar case $\hat q^{[0]}_1=0$ and $\hat q^{[0]}_2=-2$;
in the vectorial case $\hat q^{[1]}_1=\hat q^{[1]}_2=-1$;
in the pseudoscalar case $\hat q^{[12]}_1=\hat q^{[12]}_2=0$. (Theorem \ref{theo:ax} still applies.)

It is easy to see that being an axial extension automatically induces compatibility to the central axis operation:
If $\Xi$ is a vectorial Clifford conservative FQ operation, then
\begin{enumerate}
\item $\mathcal A_X\circ\Xi^{\mathrm x}=\mathcal A_X$, $\mathcal B\circ\Xi^{\mathrm x}=\mathcal B$;
\item $\mathcal A_X\circ\Xi^{\mathrm m}=\mathcal A_C$, $\mathcal B\circ\Xi^{\mathrm m}=1$;
\item $\mathcal A_X\circ\Xi^{\mathrm y}=\mathcal A_{-X}$, $\mathcal B\circ\Xi^{\mathrm y}=\mathcal B^{-1}$.
\end{enumerate}
Furthermore, one can show that $\Xi^{\mathrm x}, \Xi^{\mathrm m},\Xi^{\mathrm y}$ will automatically satisfy
scaling invariances in $\hat r_1$ and $\hat r_2$ (for all types, not only in the vectorial case).
\end{point}
\section{Examples of hyperscaling}\label{sec3}
In \cite{L2} hyperscaling had some success in describing conjugation-invariance and bivariance.
Furthermore, these  conditions  come in a great variety; hence one hopes
that they might be somewhat interesting.
On the other hand, in \cite{L2}, it was also indicated that while these conditions have a quite complex behaviour,
they are too strong to be very useful.
Here we try to enlighten this situation.
In this section we work in the formal environment, all FQ operations are meant in their formal restrictions.
\begin{point} Recall from \cite{L2} that
the FQ operation $\Xi$ satisfies the hyperscaling property
of type $(J,L,\alpha,\beta)$
in  variable  $\hat r_h$, component $[s]$, if in its expansion relative to the mixed base,
the ``decay'' identities
\begin{align}
\hat p^{[s]}_{h}&=(\alpha+\beta)\hat p^{[s]}\notag\\
\hat p^{[s]}_{h,j,\cdots}&=
\alpha \hat p^{[s]}_{j,\cdots}
-\tfrac12J\hat p^{[s]}_{6*h*j,\cdots}+(-\tfrac12-L)\hat p^{[s]}_{h*j,\cdots}\notag\\
\hat p^{[s]}_{\cdots,i,h}&=
\tfrac12J\hat p^{[s]}_{\cdots,i*h*6}+(-\tfrac12+L)\hat p^{[s]}_{\cdots,i*h}
+\beta\hat p^{[s]}_{\cdots,i}\notag\\
\hat p^{[s]}_{\cdots,i,h,j,\cdots}&=
\tfrac12J\hat p^{[s]}_{\cdots,i*h*6,j,\cdots}+(-\tfrac12+L)\hat p^{[s]}_{\cdots,i*h,j,\cdots}
-\tfrac12J\hat p^{[s]}_{\cdots,i,6*h*j,\cdots}+(-\tfrac12-L)\hat p^{[s]}_{\cdots,i,h*j,\cdots}\notag
\end{align}
hold. In what follows, this conditions will be abbreviated as $h[s](J,L,\alpha,\beta)$.

Similarly, the FQ operation $\Xi$ satisfies character degeneracy with $\pm 1$ in  variable  $\hat r_h$, component $[s]$,
if in its expansion relative to the mixed base, the identities
\[\hat p^{[s]}_{\cdots,i,\cdots}=\pm\hat p^{[s]}_{\cdots,i*h,\cdots} \]
hold.
(This is more general compared to \cite{L2}, where character degeneracy was considered only in  variable  $\hat r_6$.)
The condition above  will be abbreviated as $h[s]\langle+1\rangle$ or $h[s]\langle-1\rangle$.

It turns out, these conditions are surprisingly structured.
We start with the $L=0$ case. We will name some special types
\begin{align*}
\typI&=(1,0,1,1),&
\typII&=(1,0,1,-1),&
\typIII&=(1,0,-1,1),\\
\typIV&=(1,0,-1,-1),&
\typV&=(-1,0,0,0),&
\typVI&=(-1,0,0,1),\\
\typVII&=(-1,0,0,-1),&
\typVIII&=(1,0,1,0),&
\typIX&=(1,0,-1,0).
\end{align*}
They fit into the picture
\[\begin{array}{c|c|c|c}
L=0&\beta=1&\beta=0&\beta=-1\\\hline
\alpha=1,J=1&\typI&\typVIII&\typII\\\hline
\alpha=0,J=-1&\typVI&\typV&\typVII\\\hline
\alpha=-1,J=1&\typIII&\typIX&\typIV\\
\end{array}\quad.\]

In this terminology conjugation-invariance can be described by
\[(6[0]\typII\text{ or }6[0]\typV)\text{ and }7[0]\typII\text{ and }8[0]\typV\]
in part $[0]$;
\[(6[1]\typVI\text{ or }6[1]\typVIII)\text{ and }7[1]\typVIII\text{ and }8[1]\typVI\]
in part $[1]$;
\[(6[2]\typVI\text{ or }6[2]\typVIII)\text{ and }7[2]\typVIII\text{ and }8[2]\typVII\]
in part $[2]$;
\[(6[12]\typII\text{ or }6[12]\typV)\text{ and }7[12]\typI\text{ and }8[12]\typV\]
in part $[12]$. (We have some leverage in choosing parameters for $\hat r_6$.)
\end{point}

\begin{state}
(a) If $\Xi$ is a conjugation-invariant Clifford conservative (pseudo)scalar FQ operation satisfying the hyperscaling property
of type $(J,0,\alpha,\beta)$ in  variable  $\hat r_h$ ($h\in\{1,2,3,4,5\}$), component $[s]$ ($s\in\{0,12\}$).
Then, we claim, $\Xi$ satisfies a hyperscaling property of type
$\typI, \typII, \typIII, \typIV$ or $\typV$ in  variable  $\hat r_h$, component $[s]$.

(b) If $\Xi$ is a conjugation-invariant Clifford conservative vectorial FQ operation satisfying the hyperscaling property
of type $(J,0,\alpha,\beta)$ in  variable  $\hat r_h$ ($h\in\{1,2,3,4,5\}$), component $[s]$ ($s\in\{1,2\}$).
Then, we, claim $\Xi$ satisfies a hyperscaling property of type
$\typV, \typVI, \typVII$ or $\typIX$ in  variable  $\hat r_h$, component $[s]$.

(Remark: We do not claim that $(J,0,\alpha,\beta)$ itself is one of $\typI,\ldots,\typIX$.)
\end{state}

We call the $(5+2\times 4+5)\times 5=90$ possible hyperscaling conditions above as the principal hyperscaling conditions.
The following statement refines the situation.
It turns out that if we impose a principal hyperscaling condition $h[s]X$, then it does not only allow to eliminate $h$
from the lower indices of $\hat p^{[s]}_{\iota_1,\ldots,\iota_r}$ (beyond $6,7,8$) but implies further rules.
Ultimately, this  will allow  a reducing index set $I_{\mathrm e}\varsubsetneq\{1,\ldots,5\}\setminus\{h\}$
such that the expansion coefficients will depend on $\hat p^{[s]}_{\iota_1,\ldots,\iota_r}$ ($\iota_1,\ldots,\iota_r\in I$) which
can be  prescribed arbitrarily (with $\hat p^{[s]}=1$ in the Clifford conservative case). We term these as exact reduction sets $I_{\mathrm e}$.
When we pair hyperscaling properties with, say orthogonal invariances, we can obtain index sets $I_{\mathrm{ne}}$ such that lower
indices can be reduced to be from $I_{\mathrm{ne}}$ but subject to further conditions, i. e. free prescribability does not hold.

In the following statements we deal with Clifford conservative operations.
$\hat{\mathbf P}^{[s]}_{1\ldots5}$ means $\hat{\mathbf P}^{[s]}_1$ restricted to its first 5 entries.
\begin{state}
For scalar FQ operations  $\Xi$, the consistent constellations of principal hyperscaling properties are as follows:
\begin{enumerate}
\item $1[0]\typI \Leftrightarrow 4[0]\typI$
      $\qquad \rightsquigarrow I_{\mathrm e}=\{2,3,5\}, \hat{\mathbf P}^{[0]}_{1\ldots5}=\begin{bmatrix}2&\hat q_2&\hat q_3&2&\hat q_5\end{bmatrix}$
\item $1[0]\typII \Leftrightarrow 4[0]\typII$
      $\qquad \rightsquigarrow I_{\mathrm e}=\{2,3,5\}, \hat{\mathbf P}^{[0]}_{1\ldots5}=\begin{bmatrix}0&\hat q_2&\hat q_3&0&\hat q_5\end{bmatrix}$
\item $1[0]\typIII \Leftrightarrow 4[0]\typIII$
      $\qquad \rightsquigarrow I_{\mathrm e}=\{2,3,5\}, \hat{\mathbf P}^{[0]}_{1\ldots5}=\begin{bmatrix}0&\hat q_2&\hat q_3&0&\hat q_5\end{bmatrix}$
\item $1[0]\typIV \Leftrightarrow 4[0]\typIV$
      $\qquad \rightsquigarrow I_{\mathrm e}=\{2,3,5\}, \hat{\mathbf P}^{[0]}_{1\ldots5}=\begin{bmatrix}-2&\hat q_2&\hat q_3&-2&\hat q_5\end{bmatrix}$
\item $1[0]\typV \Leftrightarrow 3[0]\typV$
      $\qquad \rightsquigarrow I_{\mathrm e}=\{2,4,5\}, \hat{\mathbf P}^{[0]}_{1\ldots5}=\begin{bmatrix}0&\hat q_2&0&\hat q_4&\hat q_5\end{bmatrix}$
\item $2[0]\typI \Leftrightarrow 3[0]\typI$
      $\qquad \rightsquigarrow I_{\mathrm e}=\{1,4,5\}, \hat{\mathbf P}^{[0]}_{1\ldots5}=\begin{bmatrix}\hat q_1&2&2&\hat q_4&\hat q_5\end{bmatrix}$
\item $2[0]\typII \Leftrightarrow 3[0]\typII$
      $\qquad \rightsquigarrow I_{\mathrm e}=\{1,4,5\}, \hat{\mathbf P}^{[0]}_{1\ldots5}=\begin{bmatrix}\hat q_1&0&0&\hat q_4&\hat q_5\end{bmatrix}$
\item $2[0]\typIII \Leftrightarrow 3[0]\typIII$
      $\qquad \rightsquigarrow I_{\mathrm e}=\{1,4,5\}, \hat{\mathbf P}^{[0]}_{1\ldots5}=\begin{bmatrix}\hat q_1&0&0&\hat q_4&\hat q_5\end{bmatrix}$
\item $2[0]\typIV \Leftrightarrow 3[0]\typIV$
      $\qquad \rightsquigarrow I_{\mathrm e}=\{1,4,5\}, \hat{\mathbf P}^{[0]}_{1\ldots5}=\begin{bmatrix}\hat q_1&-2&-2&\hat q_4&\hat q_5\end{bmatrix}$
\item $2[0]\typV \Leftrightarrow 4[0]\typV$
      $\qquad \rightsquigarrow I_{\mathrm e}=\{1,3,5\}, \hat{\mathbf P}^{[0]}_{1\ldots5}=\begin{bmatrix}\hat q_1&0&\hat q_3&0&\hat q_5\end{bmatrix}$
\item $5[0]\typII \Leftrightarrow 1[0]\langle+1\rangle$
      $\qquad \rightsquigarrow I_{\mathrm e}=\{3,4\}, \hat{\mathbf P}^{[0]}_{1\ldots5}=\begin{bmatrix}\hat q_3&\hat q_4&\hat q_3&\hat q_4&0\end{bmatrix}$
\item $5[0]\typIII \Leftrightarrow 1[0]\langle-1\rangle$
      $\qquad \rightsquigarrow I_{\mathrm e}=\{3,4\}, \hat{\mathbf P}^{[0]}_{1\ldots5}=\begin{bmatrix}-\hat q_3&-\hat q_4&\hat q_3&\hat q_4&0\end{bmatrix}$
\item $5[0]\typV \Leftrightarrow 2[0]\langle+1\rangle$
      $\qquad \rightsquigarrow I_{\mathrm e}=\{3,4\}, \hat{\mathbf P}^{[0]}_{1\ldots5}=\begin{bmatrix}\hat q_4&\hat q_3&\hat q_3&\hat q_4&0\end{bmatrix}$
\item $1,4[0]\typII \,\&\, 2,4[0]\typV \Leftrightarrow 5[0]\langle+1\rangle$
      $\qquad \rightsquigarrow I_{\mathrm e}=\{5\}, \hat{\mathbf P}^{[0]}_{1\ldots5}=\begin{bmatrix}0&0&\hat q_5&0&\hat q_5\end{bmatrix}$
\item $1,4[0]\typIII \,\&\, 2,4[0]\typV \Leftrightarrow 5[0]\langle-1\rangle$
      $\qquad \rightsquigarrow I_{\mathrm e}=\{5\}, \hat{\mathbf P}^{[0]}_{1\ldots5}=\begin{bmatrix}0&0&-\hat q_5&0&\hat q_5\end{bmatrix}$
\item $1,4[0]\typII \,\&\, 1,3[0]\typV \Leftrightarrow 8[0]\langle+1\rangle$
      $\qquad \rightsquigarrow I_{\mathrm e}=\{5\}, \hat{\mathbf P}^{[0]}_{1\ldots5}=\begin{bmatrix}0&\hat q_5&0&0&\hat q_5\end{bmatrix}$
\item $1,4[0]\typIII \,\&\, 1,3[0]\typV \Leftrightarrow 8[0]\langle-1\rangle$
      $\qquad \rightsquigarrow I_{\mathrm e}=\{5\}, \hat{\mathbf P}^{[0]}_{1\ldots5}=\begin{bmatrix}0&-\hat q_5&0&0&\hat q_5\end{bmatrix}$
      \item $1,3[0]\typV \,\&\, 2,3[0]\typII \Leftrightarrow 6[0]\langle+1\rangle$
      $\qquad \rightsquigarrow I_{\mathrm e}=\{5\}, \hat{\mathbf P}^{[0]}_{1\ldots5}=\begin{bmatrix}0&0&0&\hat q_5&\hat q_5\end{bmatrix}$
\item $1,3[0]\typV \,\&\, 2,3[0]\typIII \Leftrightarrow 6[0]\langle-1\rangle$
      $\qquad \rightsquigarrow I_{\mathrm e}=\{5\}, \hat{\mathbf P}^{[0]}_{1\ldots5}=\begin{bmatrix}0&0&0&-\hat q_5&\hat q_5\end{bmatrix}$
\item $2,3[0]\typII \,\&\, 2,4[0]\typV \Leftrightarrow 7[0]\langle+1\rangle$
      $\qquad \rightsquigarrow I_{\mathrm e}=\{5\}, \hat{\mathbf P}^{[0]}_{1\ldots5}=\begin{bmatrix}\hat q_5&0&0&0&\hat q_5\end{bmatrix}$
\item $2,3[0]\typIII \,\&\, 2,4[0]\typV \Leftrightarrow 7[0]\langle-1\rangle$
      $\qquad \rightsquigarrow I_{\mathrm e}=\{5\}, \hat{\mathbf P}^{[0]}_{1\ldots5}=\begin{bmatrix}-\hat q_5&0&0&0&\hat q_5\end{bmatrix}$
\item $1,3[0]\typV \,\&\, 5[0]\typII$
      $\qquad \rightsquigarrow I_{\mathrm e}=\{4\}, \hat{\mathbf P}^{[0]}_{1\ldots5}=\begin{bmatrix}0&\hat q_4&0&\hat q_4&0\end{bmatrix}$
\item $1,3[0]\typV \,\&\, 5[0]\typIII$
      $\qquad \rightsquigarrow I_{\mathrm e}=\{4\}, \hat{\mathbf P}^{[0]}_{1\ldots5}=\begin{bmatrix}0&-\hat q_4&0&\hat q_4&0\end{bmatrix}$
\item $2,3[0]\typI \,\&\, 5[0]\typV$
      $\qquad \rightsquigarrow I_{\mathrm e}=\{4\}, \hat{\mathbf P}^{[0]}_{1\ldots5}=\begin{bmatrix}\hat q_4&2&2&\hat q_4&0\end{bmatrix}$
\item $2,3[0]\typII \,\&\, 5[0]\typV$
      $\qquad \rightsquigarrow I_{\mathrm e}=\{4\}, \hat{\mathbf P}^{[0]}_{1\ldots5}=\begin{bmatrix}\hat q_4&0&0&\hat q_4&0\end{bmatrix}$
\item $2,3[0]\typIII \,\&\, 5[0]\typV$
      $\qquad \rightsquigarrow I_{\mathrm e}=\{4\}, \hat{\mathbf P}^{[0]}_{1\ldots5}=\begin{bmatrix}\hat q_4&0&0&\hat q_4&0\end{bmatrix}$
\item $2,3[0]\typIV \,\&\, 5[0]\typV$
      $\qquad \rightsquigarrow I_{\mathrm e}=\{4\}, \hat{\mathbf P}^{[0]}_{1\ldots5}=\begin{bmatrix}\hat q_4&-2&-2&\hat q_4&0\end{bmatrix}$
\item $1,4[0]\typI \,\&\, 5[0]\typV$
      $\qquad \rightsquigarrow I_{\mathrm e}=\{3\}, \hat{\mathbf P}^{[0]}_{1\ldots5}=\begin{bmatrix}2&\hat q_3&\hat q_3&2&0\end{bmatrix}$
\item $1,4[0]\typII \,\&\, 5[0]\typV$
      $\qquad \rightsquigarrow I_{\mathrm e}=\{3\}, \hat{\mathbf P}^{[0]}_{1\ldots5}=\begin{bmatrix}0&\hat q_3&\hat q_3&0&0\end{bmatrix}$
\item $1,4[0]\typIII \,\&\, 5[0]\typV$
      $\qquad \rightsquigarrow I_{\mathrm e}=\{3\}, \hat{\mathbf P}^{[0]}_{1\ldots5}=\begin{bmatrix}0&\hat q_3&\hat q_3&0&0\end{bmatrix}$
\item $1,4[0]\typIV \,\&\, 5[0]\typV$
      $\qquad \rightsquigarrow I_{\mathrm e}=\{3\}, \hat{\mathbf P}^{[0]}_{1\ldots5}=\begin{bmatrix}-2&\hat q_3&\hat q_3&-2&0\end{bmatrix}$
\item $2,4[0]\typV \,\&\, 5[0]\typII$
      $\qquad \rightsquigarrow I_{\mathrm e}=\{3\}, \hat{\mathbf P}^{[0]}_{1\ldots5}=\begin{bmatrix}\hat q_3&0&\hat q_3&0&0\end{bmatrix}$
\item $2,4[0]\typV \,\&\, 5[0]\typIII$
      $\qquad \rightsquigarrow I_{\mathrm e}=\{3\}, \hat{\mathbf P}^{[0]}_{1\ldots5}=\begin{bmatrix}-\hat q_3&0&\hat q_3&0&0\end{bmatrix}$
\item $5[0]\typII \,\&\, 5[0]\typV  \Leftrightarrow 4[0]\langle+1\rangle$
      $\qquad \rightsquigarrow I_{\mathrm e}=\{3\}, \hat{\mathbf P}^{[0]}_{1\ldots5}=\begin{bmatrix}\hat q_3&\hat q_3&\hat q_3&\hat q_3&0\end{bmatrix}$
\item $5[0]\typIII \,\&\, 5[0]\typV \Leftrightarrow 4[0]\langle-1\rangle$
      $\qquad \rightsquigarrow I_{\mathrm e}=\{3\}, \hat{\mathbf P}^{[0]}_{1\ldots5}=\begin{bmatrix}-\hat q_3&\hat q_3&\hat q_3&-\hat q_3&0\end{bmatrix}$
\item $1,4[0]\typI \,\&\, 2,3[0]\typI \,\&\, 5[0]\typII$
      $\qquad \rightsquigarrow I_{\mathrm e}=\emptyset, \hat{\mathbf P}^{[0]}_{1\ldots5}=\begin{bmatrix}2&2&2&2&0\end{bmatrix}$

      and in this case $\Xi(A_1,A_2)=-A_1^2$
\item $1,4[0]\typIV \,\&\, 2,3[0]\typI\,\&\, 5[0]\typIII$
      $\qquad \rightsquigarrow I_{\mathrm e}=\emptyset, \hat{\mathbf P}^{[0]}_{1\ldots5}=\begin{bmatrix}-2&2&2&-2&0\end{bmatrix}$

      and in this case $\Xi(A_1,A_2)=-A_2^2$
\item $1,4[0]\typI \,\&\, 2,3[0]\typIV \,\&\, 5[0]\typIII$ $\qquad \rightsquigarrow I_{\mathrm e}=\emptyset, \hat{\mathbf P}^{[0]}_{1\ldots5}=\begin{bmatrix}2&-2&-2&2&0\end{bmatrix}$

      and in this case $\Xi(A_1,A_2)=-A_2^{-2}$
\item $1,4[0]\typIV \,\&\,2,3[0]\typIV \,\&\, 5[0]\typII$
      $\qquad \rightsquigarrow I_{\mathrm e}=\emptyset, \hat{\mathbf P}^{[0]}_{1\ldots5}=\begin{bmatrix}-2&-2&-2&-2&0\end{bmatrix}$

      and in this case $\Xi(A_1,A_2)=-A_1^{-2}$
\item $1,2,3,4,5[0]\typII \,\&\, 1,2,3,4,5[0]\typIII \,\&\, 1,2,3,4,5[0]\typV$
      $\quad \rightsquigarrow I_{\mathrm e}=\emptyset, \hat{\mathbf P}^{[0]}_{1\ldots5}=\begin{bmatrix}0&0&0&0&0\end{bmatrix}$

      and in this case $\Xi(A_1,A_2)=1$ .
\end{enumerate}

Remark: Conditions $5[0]\typI$,  $5[0]\typIV$, $2[0]\langle-1\rangle$, $3[0]\langle-1\rangle$ are inconsistent; $3[0]\langle+1\rangle$ is trivial.

We can  say that conditions (1--13) are the primary conditions, (14--35) are composite conditions, (36--40) are extremal conditions.
\end{state}

\begin{state}\label{san0}
For orthogonal invariant scalar FQ operations  $\Xi$, the consistent constellations of principal hyperscaling properties are as follows:
\begin{enumerate}
\item $1[0]\typV \Leftrightarrow 3[0]\typV$\label{san1}
      $\qquad \rightsquigarrow I_{\mathrm{ne}}=\{2,4,5\}, \hat{\mathbf P}^{[0]}_{1\ldots5}=\begin{bmatrix}0&\hat q_2&0&0&0\end{bmatrix}$
\item $2[0]\typI \Leftrightarrow 3[0]\typI$\label{san3}
      $\qquad \rightsquigarrow I_{\mathrm{ne}}=\{1,4,5\}, \hat{\mathbf P}^{[0]}_{1\ldots5}=\begin{bmatrix}0&2&2&0&0\end{bmatrix}$
\item $2[0]\typII \Leftrightarrow 3[0]\typII$
      $\qquad \rightsquigarrow I_{\mathrm{ne}}=\{1,4,5\}, \hat{\mathbf P}^{[0]}_{1\ldots5}=\begin{bmatrix}0&0&0&0&0\end{bmatrix}$
\item $2[0]\typIII \Leftrightarrow 3[0]\typIII$
      $\qquad \rightsquigarrow I_{\mathrm{ne}}=\{1,4,5\}, \hat{\mathbf P}^{[0]}_{1\ldots5}=\begin{bmatrix}0&0&0&0&0\end{bmatrix}$
\item $2[0]\typIV \Leftrightarrow 3[0]\typIV$\label{san4}
      $\qquad \rightsquigarrow I_{\mathrm{ne}}=\{1,4,5\}, \hat{\mathbf P}^{[0]}_{1\ldots5}=\begin{bmatrix}0&-2&-2&0&0\end{bmatrix}$
\item $2[0]\typV \Leftrightarrow 4[0]\typV$\label{san2}
      $\qquad \rightsquigarrow I_{\mathrm e}=\{3\}, \hat{\mathbf P}^{[0]}_{1\ldots5}=\begin{bmatrix}0&0&\hat q_3&0&0\end{bmatrix}$
\item $5[0]\typV \Leftrightarrow 2[0]\langle+1\rangle$
      $\qquad \rightsquigarrow I_{\mathrm e}=\{3\}, \hat{\mathbf P}^{[0]}_{1\ldots5}=\begin{bmatrix}0&\hat q_3&\hat q_3&0&0\end{bmatrix}$
\item $2,3[0]\typI \,\&\, 5[0]\typV$
      $\qquad \rightsquigarrow I_{\mathrm e}=\emptyset, \hat{\mathbf P}^{[0]}_{1\ldots5}=\begin{bmatrix}0&2&2&0&0\end{bmatrix}$

      and in this case $\Xi(A_1,A_2)=\mathcal L(A_1,A_2)\equiv-\frac12(A_1^2+A_2^2)$
\item $2,3[0]\typIV \,\&\, 5[0]\typV$
      $\qquad \rightsquigarrow I_{\mathrm e}=\emptyset, \hat{\mathbf P}^{[0]}_{1\ldots5}=\begin{bmatrix}0&-2&-2&0&0\end{bmatrix}$

      and in this case $\Xi(A_1,A_2)=\mathcal L(A_1,A_2)^{-1}\equiv-\left(\frac12(A_1^{2}+A_2^{2})\right)^{-1}$
\item $1,2,3,4,5[0]\typII \,\&\, 1,2,3,4,5[0]\typIII \,\&\, 1,2,3,4,5[0]\typV$
     $\quad \rightsquigarrow I_{\mathrm e}=\emptyset, \hat{\mathbf P}^{[0]}_{1\ldots5}=\begin{bmatrix}0&0&0&0&0\end{bmatrix}$

     and in this case $\Xi(A_1,A_2)=1$ .
\end{enumerate}

We can  say that
conditions (1--7) are the primary conditions, (8--10)  are extremal conditions.
\end{state}
\begin{state}For orthogonal invariant scalar FQ operations,

\begin{itemize}
\item an example for \ref{san0}\eqref{san1} is $\Xi(A_1,A_2)=\mathcal B(A_1,A_2)$ with values
$\hat q_2=2$, $\hat{\mathbf P}^{[0]}_{1\ldots5}=\begin{bmatrix}0&2&0&0&0\end{bmatrix}$;
\item an example for \ref{san0}\eqref{san1} is $\Xi(A_1,A_2)=\mathcal B^{-1}(A_1,A_2)$ with values
$\hat q_2=-2$, $\hat{\mathbf P}^{[0]}_{1\ldots5}=\begin{bmatrix}0&-2&0&0&0\end{bmatrix}$;
\item examples for \ref{san0}\eqref{san3} are $\Xi=\mathcal L^{\mathrm{fx}},\mathcal L^{\mathrm{x}},\mathcal L^{\mathrm{cx}}$;
\item examples for \ref{san0}\eqref{san4} are $\Xi=(\mathcal L^{\mathrm{fx}})^{-1},(\mathcal L^{\mathrm{x}})^{-1},(\mathcal L^{\mathrm{cx}})^{-1}$;
\item an  example for \ref{san0}\eqref{san2} is $\Xi(A_1,A_2)=|\mathcal D(A_1,A_2)|$ with values $\hat q_3=2$,
$\hat{\mathbf P}^{[0]}_{1\ldots5}=\begin{bmatrix}0&0&2&0&0\end{bmatrix}$;
\item an  example for \ref{san0}\eqref{san2} is $\Xi(A_1,A_2)=|\mathcal D(A_1,A_2)|^{-1}$ with  values $\hat q_3=-2$,
$\hat{\mathbf P}^{[0]}_{1\ldots5}=\begin{bmatrix}0&0&-2&0&0\end{bmatrix}$.

\end{itemize}
\end{state}

\begin{state}
For vectorial FQ operations $\Xi$, which can be assumed to be symmetric,  the consistent constellations of principal hyperscaling properties in component $[1]$ are as follows:
\begin{enumerate}
\item $1[1]\typVI \Leftrightarrow 3[1]\typVI$
      $\qquad \rightsquigarrow I_{\mathrm e}=\{2,4,5\}, \hat{\mathbf P}^{[1]}_{1\ldots5}=\begin{bmatrix}1&\hat q_2&1&\hat q_4&\hat q_5\end{bmatrix}$
\item $1[1]\typVII \Leftrightarrow 3[1]\typVII$
      $\qquad \rightsquigarrow I_{\mathrm e}=\{2,4,5\}, \hat{\mathbf P}^{[1]}_{1\ldots5}=\begin{bmatrix}-1&\hat q_2&-1&\hat q_4&\hat q_5\end{bmatrix}$
\item $1[1]\typVIII \Leftrightarrow 4[1]\typVIII$
      $\qquad \rightsquigarrow I_{\mathrm e}=\{2,3,5\}, \hat{\mathbf P}^{[1]}_{1\ldots5}=\begin{bmatrix}1&\hat q_2&\hat q_3&1&\hat q_5\end{bmatrix}$
\item $1[1]\typIX \Leftrightarrow 4[1]\typIX$
      $\qquad \rightsquigarrow I_{\mathrm e}=\{2,3,5\}, \hat{\mathbf P}^{[1]}_{1\ldots5}=\begin{bmatrix}-1&\hat q_2&\hat q_3&-1&\hat q_5\end{bmatrix}$
\item $2[1]\typVI \Leftrightarrow 4[1]\typVI$
      $\qquad \rightsquigarrow I_{\mathrm e}=\{1,3,5\}, \hat{\mathbf P}^{[1]}_{1\ldots5}=\begin{bmatrix}\hat q_1&1&\hat q_3&1&\hat q_5\end{bmatrix}$
\item $2[1]\typVII \Leftrightarrow 4[1]\typVII$
      $\qquad \rightsquigarrow I_{\mathrm e}=\{1,3,5\}, \hat{\mathbf P}^{[1]}_{1\ldots5}=\begin{bmatrix}\hat q_1&-1&\hat q_3&-1&\hat q_5\end{bmatrix}$
\item $2[1]\typVIII \Leftrightarrow 3[1]\typVIII$
      $\qquad \rightsquigarrow I_{\mathrm e}=\{1,4,5\}, \hat{\mathbf P}^{[1]}_{1\ldots5}=\begin{bmatrix}\hat q_1&1&1&\hat q_4&\hat q_5\end{bmatrix}$
\item $2[1]\typIX \Leftrightarrow 3[1]\typIX$
      $\qquad \rightsquigarrow I_{\mathrm e}=\{1,4,5\}, \hat{\mathbf P}^{[1]}_{1\ldots5}=\begin{bmatrix}\hat q_1&-1&-1&\hat q_4&\hat q_5\end{bmatrix}$
\item $5[1]\typVI \Leftrightarrow 2[1]\langle+1\rangle$
      $\qquad \rightsquigarrow I_{\mathrm e}=\{3,4\}, \hat{\mathbf P}^{[1]}_{1\ldots5}=\begin{bmatrix}\hat q_4&\hat q_3&\hat q_3&\hat q_4&1\end{bmatrix}$
\item $5[1]\typVIII \Leftrightarrow 1[1]\langle+1\rangle$
      $\qquad \rightsquigarrow I_{\mathrm e}=\{3,4\}, \hat{\mathbf P}^{[1]}_{1\ldots5}=\begin{bmatrix}\hat q_3&\hat q_4&\hat q_3&\hat q_4&1\end{bmatrix}$
 \item $1,3[1]\typVI \,\&\, 2,3[1]\typVIII\Leftrightarrow6[1]\langle+1\rangle$
      $\qquad \rightsquigarrow I_{\mathrm e}=\{5\}, \hat{\mathbf P}^{[1]}_{1\ldots5}=\begin{bmatrix}1&1&1&\hat q_5&\hat q_5\end{bmatrix}$
\item $1,3[1]\typVII \,\&\, 2,3[1]\typIX\Leftrightarrow6[1]\langle-1\rangle$
      $\qquad \rightsquigarrow I_{\mathrm e}=\{5\}, \hat{\mathbf P}^{[1]}_{1\ldots5}=\begin{bmatrix}-1&-1&-1&-\hat q_5&\hat q_5\end{bmatrix}$
\item $2,4[1]\typVI \,\&\, 2,3[1]\typVIII\Leftrightarrow7[1]\langle+1\rangle$
      $\qquad \rightsquigarrow I_{\mathrm e}=\{5\}, \hat{\mathbf P}^{[1]}_{1\ldots5}=\begin{bmatrix}\hat q_5&1&1&1&\hat q_5\end{bmatrix}$
\item $2,4[1]\typVII \,\&\, 2,3[1]\typIX\Leftrightarrow7[1]\langle-1\rangle$
      $\qquad \rightsquigarrow I_{\mathrm e}=\{5\}, \hat{\mathbf P}^{[1]}_{1\ldots5}=\begin{bmatrix}-\hat q_5&-1&-1&-1&\hat q_5\end{bmatrix}$
\item $1,3[1]\typVI \,\&\, 1,4[1]\typVIII\Leftrightarrow8[1]\langle+1\rangle$
      $\qquad \rightsquigarrow I_{\mathrm e}=\{5\}, \hat{\mathbf P}^{[1]}_{1\ldots5}=\begin{bmatrix}1&\hat q_5&1&1&\hat q_5\end{bmatrix}$
\item $1,3[1]\typVII \,\&\, 1,4[1]\typIX\Leftrightarrow8[1]\langle-1\rangle$
      $\qquad \rightsquigarrow I_{\mathrm e}=\{5\}, \hat{\mathbf P}^{[1]}_{1\ldots5}=\begin{bmatrix}-1&-\hat q_5&-1&-1&\hat q_5\end{bmatrix}$
\item $1,4[1]\typVIII \,\&\, 2,4[1]\typVI\Leftrightarrow5[1]\langle+1\rangle$
      $\qquad \rightsquigarrow I_{\mathrm e}=\{5\}, \hat{\mathbf P}^{[1]}_{1\ldots5}=\begin{bmatrix}1&1&\hat q_5&1&\hat q_5\end{bmatrix}$
\item $1,4[1]\typIX \,\&\, 2,4[1]\typVII\Leftrightarrow5[1]\langle-1\rangle$
      $\qquad \rightsquigarrow I_{\mathrm e}=\{5\}, \hat{\mathbf P}^{[1]}_{1\ldots5}=\begin{bmatrix}-1&-1&-\hat q_5&-1&\hat q_5\end{bmatrix}$
\item $1,3[1]\typVI \,\&\, 5[1]\typVIII$
      $\qquad \rightsquigarrow I_{\mathrm e}=\{4\}, \hat{\mathbf P}^{[1]}_{1\ldots5}=\begin{bmatrix}1&\hat q_4&1&\hat q_4&1\end{bmatrix}$
\item $1,3[1]\typVII \,\&\, 5[1]\typVIII$
      $\qquad \rightsquigarrow I_{\mathrm e}=\{4\}, \hat{\mathbf P}^{[1]}_{1\ldots5}=\begin{bmatrix}-1&\hat q_4&-1&\hat q_4&1\end{bmatrix}$
\item $2,3[1]\typVIII \,\&\, 5[1]\typVI$
      $\qquad \rightsquigarrow I_{\mathrm e}=\{4\}, \hat{\mathbf P}^{[1]}_{1\ldots5}=\begin{bmatrix}\hat q_4&1&1&\hat q_4&1\end{bmatrix}$
\item $2,3[1]\typIX \,\&\, 5[1]\typVI$
      $\qquad \rightsquigarrow I_{\mathrm e}=\{4\}, \hat{\mathbf P}^{[1]}_{1\ldots5}=\begin{bmatrix}\hat q_4&-1&-1&\hat q_4&1\end{bmatrix}$
\item $1,4[1]\typVIII \,\&\, 5[1]\typVI$
      $\qquad \rightsquigarrow I_{\mathrm e}=\{3\}, \hat{\mathbf P}^{[1]}_{1\ldots5}=\begin{bmatrix}1&\hat q_3&\hat q_3&1&1\end{bmatrix}$
\item $1,4[1]\typIX \,\&\, 5[1]\typVI$
      $\qquad \rightsquigarrow I_{\mathrm e}=\{3\}, \hat{\mathbf P}^{[1]}_{1\ldots5}=\begin{bmatrix}-1&\hat q_3&\hat q_3&-1&1\end{bmatrix}$
\item $2,4[1]\typVI \,\&\, 5[1]\typVIII$
      $\qquad \rightsquigarrow I_{\mathrm e}=\{3\}, \hat{\mathbf P}^{[1]}_{1\ldots5}=\begin{bmatrix}\hat q_3&1&\hat q_3&1&1\end{bmatrix}$
\item $2,4[1]\typVII \,\&\, 5[1]\typVIII$
      $\qquad \rightsquigarrow I_{\mathrm e}=\{3\}, \hat{\mathbf P}^{[1]}_{1\ldots5}=\begin{bmatrix}-\hat q_3&-1&-\hat q_3&-1&1\end{bmatrix}$
\item $5[1]\typVI \,\&\, 5[1]\typVIII\Leftrightarrow4[1]\langle+1\rangle$
      $\qquad \rightsquigarrow I_{\mathrm e}=\{3\}, \hat{\mathbf P}^{[1]}_{1\ldots5}=\begin{bmatrix}\hat q_3&\hat q_3&\hat q_3&\hat q_3&1\end{bmatrix}$
\item $1,2,3,4,5[1]\typVI \,\&\, 1,2,3,4,5[1]\typVIII$
      $\qquad \rightsquigarrow I_{\mathrm e}=\emptyset, \hat{\mathbf P}^{[1]}_{1\ldots5}=\begin{bmatrix}1&1&1&1&1\end{bmatrix}$

      and in this case $\Xi(A_1,A_2)=\Id(A_1,A_2)=(A_1,A_2)$

\item $1,2,3,4[1]\typVII \,\&\, 1,2,3,4[1]\typIX\,\&\,5[1]\typVI\,\&\,5[1]\typVIII$
      $\quad \rightsquigarrow I_{\mathrm e}=\emptyset, $
      $\hat{\mathbf P}^{[1]}_{1\ldots5}=\begin{bmatrix}-1&\!-1&\!-1&\!-1&1\end{bmatrix}$

      and in this case $\Xi(A_1,A_2)=(-A_1^{-1},-A_2^{-1})$
\end{enumerate}

Remark: Conditions $5[1]\typVII$,  $5[1]\typIX$, $1[1]\langle-1\rangle$, $2[1]\langle-1\rangle$, $3[1]\langle-1\rangle$, $4[1]\langle-1\rangle$  are inconsistent; $3[1]\langle+1\rangle$ is trivial.

We can say that
conditions (1--10) are the primary conditions, (11--27) are composite conditions, (28--29) are extremal conditions.
\end{state}

\begin{state}\label{van0}
For orthogonal-invariant vectorial FQ operations $\Xi$ (which are necessarily symmetric),  the consistent constellations of principal hyperscaling properties in component $[1]$ are as follows:
\begin{enumerate}
\item $1[1]\typVI \Leftrightarrow 3[1]\typVI$ \label{van2}
      $\qquad \rightsquigarrow I_{\mathrm{ne}}=\{2,4,5\}, \hat{\mathbf P}^{[1]}_{1\ldots5}=\begin{bmatrix}1&\hat q_2&1&\hat q_5&\hat q_5\end{bmatrix}$
\item $1[1]\typVII \Leftrightarrow 3[1]\typVII$\label{van6}
      $\qquad \rightsquigarrow I_{\mathrm{ne}}=\{2,4,5\}, \hat{\mathbf P}^{[1]}_{1\ldots5}=\begin{bmatrix}-1&\hat q_2&-1&\hat q_5&\hat q_5\end{bmatrix}$
\item $2[1]\typVIII \Leftrightarrow 3[1]\typVIII$ \label{van1}
      $\qquad \rightsquigarrow I_{\mathrm{ne}}=\{1,4,5\}, \hat{\mathbf P}^{[1]}_{1\ldots5}=\begin{bmatrix}\hat q_1&1&1&\hat q_5&\hat q_5\end{bmatrix}$
\item $2[1]\typIX \Leftrightarrow 3[1]\typIX$\label{van7}
      $\qquad \rightsquigarrow I_{\mathrm{ne}}=\{1,4,5\}, \hat{\mathbf P}^{[1]}_{1\ldots5}=\begin{bmatrix}\hat q_1&-1&-1&\hat q_5&\hat q_5\end{bmatrix}$
\item $2[1]\typVI \Leftrightarrow 4[1]\typVI$\label{van4}
      $\qquad \rightsquigarrow I_{\mathrm e}=\{3\}, \hat{\mathbf P}^{[1]}_{1\ldots5}=\begin{bmatrix}1&1&\hat q_3&1&1\end{bmatrix}$
\item $2[1]\typVII \Leftrightarrow 4[1]\typVII$
      $\qquad \rightsquigarrow I_{\mathrm e}=\{3\}, \hat{\mathbf P}^{[1]}_{1\ldots5}=\begin{bmatrix}1&-1&\hat q_3&-1&-1\end{bmatrix}$
\item $1[1]\typVIII \Leftrightarrow 4[1]\typVIII$\label{van3}
      $\qquad \rightsquigarrow I_{\mathrm e}=\{3\}, \hat{\mathbf P}^{[1]}_{1\ldots5}=\begin{bmatrix}1&1&\hat q_3&1&1\end{bmatrix}$
\item $1[1]\typIX \Leftrightarrow 4[1]\typIX$
      $\qquad \rightsquigarrow I_{\mathrm e}=\{3\}, \hat{\mathbf P}^{[1]}_{1\ldots5}=\begin{bmatrix}-1&1&\hat q_3&-1&-1\end{bmatrix}$
\item $5[1]\typVI \Leftrightarrow 2[1]\langle+1\rangle$
      $\qquad \rightsquigarrow I_{\mathrm e}=\{3\}, \hat{\mathbf P}^{[1]}_{1\ldots5}=\begin{bmatrix}1&\hat q_3&\hat q_3&1&1\end{bmatrix}$
\item $5[1]\typVIII \Leftrightarrow 1[1]\langle+1\rangle$
      $\qquad \rightsquigarrow I_{\mathrm e}=\{3\}, \hat{\mathbf P}^{[1]}_{1\ldots5}=\begin{bmatrix}\hat q_3&1&\hat q_3&1&1\end{bmatrix}$
\item $1,3[1]\typVI \,\&\, 2,3[1]\typVIII\Leftrightarrow6[1]\langle+1\rangle$
      $\qquad \rightsquigarrow I_{\mathrm{ne}}=\{5\}, \hat{\mathbf P}^{[1]}_{1\ldots5}=\begin{bmatrix}1&1&1&\hat q_5&\hat q_5\end{bmatrix}$

      and in this case $\Xi(A_1,A_2)=\hat q_5\cdot(A_1,A_2)+(1-\hat q_5)\cdot\mathcal O^{\mathrm{fSy}}(A_1,A_2)$
\item $1,3[1]\typVII \,\&\, 2,3[1]\typIX\Leftrightarrow6[1]\langle-1\rangle$
      $\qquad \rightsquigarrow I_{\mathrm e}=\emptyset, \hat{\mathbf P}^{[1]}_{1\ldots5}=\begin{bmatrix}-1&-1&-1&0&0\end{bmatrix}$

      and in this case $\Xi(A_1,A_2)={\mathcal O}^{\mathrm{afSy}}(A_1,A_2)$
\item $2,4[1]\typVII \,\&\, 2,3[1]\typIX\Leftrightarrow7[1]\langle-1\rangle$
      $\qquad \rightsquigarrow I_{\mathrm e}=\emptyset, \hat{\mathbf P}^{[1]}_{1\ldots5}=\begin{bmatrix}1&-1&-1&-1&-1\end{bmatrix}$

      and in this case $\Xi(A_1,A_2)=\mathcal D(A_1,A_2)^{-1}\cdot(A_2,-A_1)$
\item $1,3[1]\typVII \,\&\, 1,4[1]\typIX\Leftrightarrow8[1]\langle-1\rangle$
      $\qquad \rightsquigarrow I_{\mathrm e}=\emptyset, \hat{\mathbf P}^{[1]}_{1\ldots5}=\begin{bmatrix}-1&1&-1&-1&-1\end{bmatrix}$

      and in this case $\Xi(A_1,A_2)=(-A_2,A_1)\cdot\mathcal D(A_1,A_2)^{-1}$
\item $1,3[1]\typVII \,\&\, 5[1]\typVIII$
      $\qquad \rightsquigarrow I_{\mathrm e}=\emptyset, \hat{\mathbf P}^{[1]}_{1\ldots5}=\begin{bmatrix}-1&1&-1&1&1\end{bmatrix}$

      and in this case $\Xi(A_1,A_2)=(A_1,A_2)\cdot\mathcal L(A_1,A_2)^{-1}$
      \item $2,3[1]\typIX \,\&\, 5[1]\typVI$
      $\qquad \rightsquigarrow I_{\mathrm e}=\emptyset, \hat{\mathbf P}^{[1]}_{1\ldots5}=\begin{bmatrix}1&-1&-1&1&1\end{bmatrix}$

      and in this case $\Xi(A_1,A_2)=\mathcal L(A_1,A_2)^{-1}\cdot(A_1,A_2)$
\item $1,2,3,4,5[1]\typVI \,\&\, 1,2,3,4,5[1]\typVIII$
      $\qquad \rightsquigarrow I_{\mathrm e}=\emptyset, \hat{\mathbf P}^{[1]}_{1\ldots5}=\begin{bmatrix}1&1&1&1&1\end{bmatrix}$

      and in this case $\Xi(A_1,A_2)=\Id(A_1,A_2)=(A_1,A_2)$
\end{enumerate}

Moreover, if $\Xi$ is also transposition invariant, then it is  $\hat q_5\cdot\Id+(1-\hat q_5)\cdot\mathcal O^{\mathrm{fSy}}$
or $\mathcal O^{\mathrm{afSy}}$.

We can  say that
conditions (1--10) are the primary conditions, (11) is a composite condition, (12--17) are extremal conditions.
\end{state}
\begin{state}For orthogonal invariant vectorial FQ operations,
\begin{itemize}

\item an example for \ref{van0}\eqref{van2} is $\Xi(A_1,A_2)=\mathcal M_R(A_1,A_2)\equiv\mathcal B(A_1,A_2)^{-1}\cdot(A_1,A_2)$ with
$\hat q_2=-1$, $\hat{\mathbf P}^{[1]}_{1\ldots5}=\begin{bmatrix}1&-1&1&1&1\end{bmatrix}$;
\item an example for \ref{van0}\eqref{van1} is $\Xi(A_1,A_2)=\mathcal M_L(A_1,A_2)\equiv(A_1,A_2)\cdot\mathcal B(A_1,A_2)^{-1}$ with
$\hat q_1=-1$, $\hat{\mathbf P}^{[1]}_{1\ldots5}=\begin{bmatrix}-1&1&1&1&1\end{bmatrix}$;

\item an  example for \ref{van0}\eqref{van2} is $\Xi(A_1,A_2)=\mathcal T_ {RR}(A_1,A_2):=\mathcal A_R(A_1,A_2)\cdot(-A_2,A_1)$ with
$\hat q_2=-3,\hat q_5=-1$, $ \hat{\mathbf P}^{[1]}_{1\ldots5}=\begin{bmatrix}1&-3&1&-1&-1\end{bmatrix}$;
\item an example for \ref{van0}\eqref{van1} is $\Xi(A_1,A_2)=\mathcal T_ {LL}(A_1,A_2):=(A_2,-A_1)\cdot\mathcal A_L(A_1,A_2)$ with
$\hat q_1=-3,\hat q_5=-1$, $\hat{\mathbf P}^{[1]}_{1\ldots5}=\begin{bmatrix}-3&1&1&-1&-1\end{bmatrix}$;

\item an  example for \ref{van0}\eqref{van2} is $\Xi(A_1,A_2)=\mathcal T_ {R}(A_1,A_2):=\mathcal A_C(A_1,A_2)\cdot(-A_2,A_1)$ with
$\hat q_2=-1,\hat q_5=-1$, $ \hat{\mathbf P}^{[1]}_{1\ldots5}=\begin{bmatrix}1&-1&1&-1&-1\end{bmatrix}$;
\item an example for \ref{van0}\eqref{van1} is $\Xi(A_1,A_2)=\mathcal T_ {L}(A_1,A_2):=(A_2,-A_1)\cdot\mathcal A_C(A_1,A_2)$ with
$\hat q_1=-1,\hat q_5=-1$, $\hat{\mathbf P}^{[1]}_{1\ldots5}=\begin{bmatrix}-1&1&1&-1&-1\end{bmatrix}$;

\item examples for \ref{van0}\eqref{van6} are
$\Xi(A_1,A_2)=(A_1,A_2)\cdot\mathcal L^{\mathrm{fx}}(A_1,A_2)^{-1},$

$\Xi(A_1,A_2)=(A_1,A_2)\cdot\mathcal L^{\mathrm{x}}(A_1,A_2)^{-1},$
$\Xi(A_1,A_2)=(A_1,A_2)\cdot\mathcal L^{\mathrm{cx}}(A_1,A_2)^{-1}$ with

$\hat q_2=1,\hat q_5=1$, $\hat{\mathbf P}^{[1]}_{1\ldots5}=\begin{bmatrix}-1&1&-1&1&1\end{bmatrix}$;

\item examples for \ref{van0}\eqref{van7} are
 $\Xi(A_1,A_2)=\mathcal L^{\mathrm{fx}}(A_1,A_2)^{-1}\cdot(A_1,A_2),$

 $\Xi(A_1,A_2)=\mathcal L^{\mathrm{x}}(A_1,A_2)^{-1}\cdot(A_1,A_2),$
 $\Xi(A_1,A_2)=\mathcal L^{\mathrm{cx}}(A_1,A_2)^{-1}\cdot(A_1,A_2)$ with

$\hat q_1=1,\hat q_5=1$, $\hat{\mathbf P}^{[1]}_{1\ldots5}=\begin{bmatrix}1&-1&-1&1&1\end{bmatrix}$;

\item examples for \ref{van0}\eqref{van6} are $\Xi=\mathcal I_L^{\mathrm{fy}},\mathcal I_L^{\mathrm{y}},\mathcal I_L^{\mathrm{cy}} $

with
$\hat q_2=-1,\hat q_5=-1$, $\hat{\mathbf P}^{[1]}_{1\ldots5}=\begin{bmatrix}-1&-1&-1&-1&-1\end{bmatrix}$;
\item examples for \ref{van0}\eqref{van7} are $\Xi=\mathcal I_R^{\mathrm{fy}},\mathcal I_R^{\mathrm{y}},\mathcal I_R^{\mathrm{cy}} $

with
$\hat q_1=-1,\hat q_5=-1$, $\hat{\mathbf P}^{[1]}_{1\ldots5}=\begin{bmatrix}-1&-1&-1&-1&-1\end{bmatrix}$;

\item examples for \ref{van0}\eqref{van6} are
$\Xi(A_1,A_2)=(-A_2 ,A_1)\cdot\mathcal D^{\mathrm{fx}}(A_1,A_2)^{-1},$

$\Xi(A_1,A_2)=(-A_2,A_1)\cdot\mathcal D^{\mathrm{x}}(A_1,A_2)^{-1},$
$\Xi(A_1,A_2)=(-A_2,A_1)\cdot\mathcal D^{\mathrm{cx}}(A_1,A_2)^{-1}$

with $\hat q_2=1,\hat q_5=-1$, $\hat{\mathbf P}^{[1]}_{1\ldots5}=\begin{bmatrix}-1&1&-1&-1&-1\end{bmatrix}$;

\item examples for \ref{van0}\eqref{van7} are
$\Xi(A_1,A_2)=\mathcal D^{\mathrm{fx}}(A_1,A_2)^{-1}\cdot(A_2,-A_1),$

$\Xi(A_1,A_2)=\mathcal D^{\mathrm{x}}(A_1,A_2)^{-1}\cdot(A_2,-A_1),$
$\Xi(A_1,A_2)=\mathcal D^{\mathrm{cx}}(A_1,A_2)^{-1}\cdot(A_2,-A_1)$

with $\hat q_1=1,\hat q_5=-1$, $\hat{\mathbf P}^{[1]}_{1\ldots5}=\begin{bmatrix}1&-1&-1&-1&-1\end{bmatrix}$;

\item an example for \ref{van0}\eqref{van4} is $\Xi(A_1,A_2)=\mathcal F_R(A_1,A_2):=|\mathcal D(A_1,A_2)|^{-1}\cdot(A_1,A_2)$ with
$\hat q_3=-1$, $\hat{\mathbf P}^{[1]}_{1\ldots5}=\begin{bmatrix}1&1&-1&1&1\end{bmatrix}$;
\item an example for \ref{van0}\eqref{van3} is $\Xi(A_1,A_2)=\mathcal F_L(A_1,A_2):=(A_1,A_2)\cdot|\mathcal D(A_1,A_2)|^{-1}$ with
$\hat q_3=-1$, $\hat{\mathbf P}^{[1]}_{1\ldots5}=\begin{bmatrix}1&1&-1&1&1\end{bmatrix}$.
\end{itemize}
\end{state}

\begin{state}
For pseudoscalar FQ operations  $\Xi$, the consistent constellations of principal hyperscaling properties are as follows:
\begin{enumerate}
\item $1[12]\typI \Leftrightarrow 4[12]\typII$
      $\qquad \rightsquigarrow I_{\mathrm e}=\{2,3,5\}, \hat{\mathbf P}^{[12]}_{1\ldots5}=\begin{bmatrix}2&\hat q_2&\hat q_3&0&\hat q_5\end{bmatrix}$
\item $1[12]\typII \Leftrightarrow 4[12]\typI$
      $\qquad \rightsquigarrow I_{\mathrm e}=\{2,3,5\}, \hat{\mathbf P}^{[12]}_{1\ldots5}=\begin{bmatrix}0&\hat q_2&\hat q_3&2&\hat q_5\end{bmatrix}$
\item $1[12]\typIII \Leftrightarrow 4[12]\typIV$
      $\qquad \rightsquigarrow I_{\mathrm e}=\{2,3,5\}, \hat{\mathbf P}^{[12]}_{1\ldots5}=\begin{bmatrix}0&\hat q_2&\hat q_3&-2&\hat q_5\end{bmatrix}$
\item $1[12]\typIV \Leftrightarrow 4[12]\typIII$
      $\qquad \rightsquigarrow I_{\mathrm e}=\{2,3,5\}, \hat{\mathbf P}^{[12]}_{1\ldots5}=\begin{bmatrix}-2&\hat q_2&\hat q_3&0&\hat q_5\end{bmatrix}$
\item $1[12]\typV \Leftrightarrow 3[12]\typV$
      $\qquad \rightsquigarrow I_{\mathrm e}=\{2,4,5\}, \hat{\mathbf P}^{[12]}_{1\ldots5}=\begin{bmatrix}0&\hat q_2&0&\hat q_4&\hat q_5\end{bmatrix}$
\item $2[12]\typI \Leftrightarrow 3[12]\typII$
      $\qquad \rightsquigarrow I_{\mathrm e}=\{1,4,5\}, \hat{\mathbf P}^{[12]}_{1\ldots5}=\begin{bmatrix}\hat q_1&2&0&\hat q_4&\hat q_5\end{bmatrix}$
\item $2[12]\typII \Leftrightarrow 3[12]\typI$
      $\qquad \rightsquigarrow I_{\mathrm e}=\{1,4,5\}, \hat{\mathbf P}^{[12]}_{1\ldots5}=\begin{bmatrix}\hat q_1&0&2&\hat q_4&\hat q_5\end{bmatrix}$
\item $2[12]\typIII \Leftrightarrow 3[12]\typIV$
      $\qquad \rightsquigarrow I_{\mathrm e}=\{1,4,5\}, \hat{\mathbf P}^{[12]}_{1\ldots5}=\begin{bmatrix}\hat q_1&0&-2&\hat q_4&\hat q_5\end{bmatrix}$
\item $2[12]\typIV \Leftrightarrow 3[12]\typIII$
      $\qquad \rightsquigarrow I_{\mathrm e}=\{1,4,5\}, \hat{\mathbf P}^{[12]}_{1\ldots5}=\begin{bmatrix}\hat q_1&-2&0&\hat q_4&\hat q_5\end{bmatrix}$
\item $2[12]\typV\Leftrightarrow 4[12]\typV$
      $\qquad \rightsquigarrow I_{\mathrm e}=\{1,3,5\}, \hat{\mathbf P}^{[12]}_{1\ldots5}=\begin{bmatrix}\hat q_1&0&\hat q_3&0&\hat q_5\end{bmatrix}$
\item $5[12]\typI\Leftrightarrow 1[12]\langle+1\rangle$
      $\qquad \rightsquigarrow I_{\mathrm e}=\{3,4\}, \hat{\mathbf P}^{[12]}_{1\ldots5}=\begin{bmatrix}\hat q_3&\hat q_4&\hat q_3&\hat q_4&0\end{bmatrix}$
 \item $5[12]\typIV\Leftrightarrow 1[12]\langle-1\rangle$
      $\qquad \rightsquigarrow I_{\mathrm e}=\{3,4\}, \hat{\mathbf P}^{[12]}_{1\ldots5}=\begin{bmatrix}-\hat q_3&-\hat q_4&\hat q_3&\hat q_4&0\end{bmatrix}$

\item $1[12]\typI,4[12]\typII \,\&\,2[12]\typV,4[12]\typV\Leftrightarrow5[12]\langle+1\rangle$
      $\, \rightsquigarrow I_{\mathrm e}=\{5\}, \hat{\mathbf P}^{[12]}_{1\ldots5}=\begin{bmatrix}2&0&\hat q_5&0&\hat q_5\end{bmatrix}$
\item $1[12]\typIV,4[12]\typIII \,\&\,2[12]\typV,4[12]\typV\Leftrightarrow5[12]\langle-1\rangle$
      $\, \rightsquigarrow I_{\mathrm e}=\{5\}, \hat{\mathbf P}^{[12]}_{1\ldots5}=\begin{bmatrix}-2&0&-\hat q_5&0&\hat q_5\end{bmatrix}$
\item $1[12]\typV,3[12]\typV \,\&\,2[12]\typI,3[12]\typII\Leftrightarrow6[12]\langle+1\rangle$
      $\, \rightsquigarrow I_{\mathrm e}=\{5\}, \hat{\mathbf P}^{[12]}_{1\ldots5}=\begin{bmatrix}0&2&0&\hat q_5&\hat q_5\end{bmatrix}$
\item $1[12]\typV,3[12]\typV \,\&\,2[12]\typIV,3[12]\typIII\Leftrightarrow6[12]\langle-1\rangle$
      $\, \rightsquigarrow I_{\mathrm e}=\{5\}, \hat{\mathbf P}^{[12]}_{1\ldots5}=\begin{bmatrix}0&-2&0&-\hat q_5&\hat q_5\end{bmatrix}$
\item $2[12]\typII,3[12]\typI \,\&\,2[12]\typV,4[12]\typV\Leftrightarrow7[12]\langle+1\rangle$
      $\, \rightsquigarrow I_{\mathrm e}=\{5\}, \hat{\mathbf P}^{[12]}_{1\ldots5}=\begin{bmatrix}\hat q_5&0&2&0&\hat q_5\end{bmatrix}$
\item $2[12]\typIII,3[12]\typIV \,\&\,2[12]\typV,4[12]\typV\Leftrightarrow7[12]\langle-1\rangle$
      $\, \rightsquigarrow I_{\mathrm e}=\{5\}, \hat{\mathbf P}^{[12]}_{1\ldots5}=\begin{bmatrix}-\hat q_5&0&-2&0&\hat q_5\end{bmatrix}$
\item $1[12]\typII,4[12]\typI \,\&\,1[12]\typV,3[12]\typV\Leftrightarrow8[12]\langle+1\rangle$
      $\, \rightsquigarrow I_{\mathrm e}=\{5\}, \hat{\mathbf P}^{[12]}_{1\ldots5}=\begin{bmatrix}0&\hat q_5&0&2&\hat q_5\end{bmatrix}$
\item $1[12]\typIII,4[12]\typIV \,\&\,1[12]\typV,3[12]\typV\Leftrightarrow8[12]\langle-1\rangle$
      $\, \rightsquigarrow I_{\mathrm e}=\{5\}, \hat{\mathbf P}^{[12]}_{1\ldots5}=\begin{bmatrix}0&-\hat q_5&0&-2&\hat q_5\end{bmatrix}$
\item $1[12]\typV,3[12]\typV \,\&\,5[12]\typI$
      $\qquad \rightsquigarrow I_{\mathrm e}=\{4\}, \hat{\mathbf P}^{[12]}_{1\ldots5}=\begin{bmatrix}0&\hat q_4&0&\hat q_4&2\end{bmatrix}$
\item $1[12]\typV,3[12]\typV \,\&\,5[12]\typIV$
      $\qquad \rightsquigarrow I_{\mathrm e}=\{4\}, \hat{\mathbf P}^{[12]}_{1\ldots5}=\begin{bmatrix}0&-\hat q_4&0&\hat q_4&-2\end{bmatrix}$

\item $2[12]\typV,4[12]\typV \,\&\,5[12]\typI$
      $\qquad \rightsquigarrow I_{\mathrm e}=\{3\}, \hat{\mathbf P}^{[12]}_{1\ldots5}=\begin{bmatrix}\hat q_3&0&\hat q_3&0&2\end{bmatrix}$
\item $2[12]\typV,4[12]\typV \,\&\,5[12]\typIV$
      $\qquad \rightsquigarrow I_{\mathrm e}=\{3\}, \hat{\mathbf P}^{[12]}_{1\ldots5}=\begin{bmatrix}-\hat q_3&0&\hat q_3&0&-2\end{bmatrix}$

\item $1[12]\typI,4[12]\typII \,\&\,2[12]\typII,3[12]\typI\,\&\,5[12]\typI$
      $\qquad \rightsquigarrow I_{\mathrm e}=\emptyset, \hat{\mathbf P}^{[12]}_{1\ldots5}=\begin{bmatrix}2&0&2&0&2\end{bmatrix}$

      and in this case $\Xi(A_1,A_2)=A_1 A_2$
\item $1[12]\typI,4[12]\typII \,\&\,2[12]\typIII,3[12]\typIV\,\&\,5[12]\typIV$
      $\qquad \rightsquigarrow I_{\mathrm e}=\emptyset, \hat{\mathbf P}^{[12]}_{1\ldots5}=\begin{bmatrix}2&0&-2&0&-2\end{bmatrix}$

      and in this case $\Xi(A_1,A_2)=-A_2^{-1}A_1^{-1}$
\item $1[12]\typIV,4[12]\typIII \,\&\,2[12]\typII,3[12]\typI\,\&\,5[12]\typIV$
      $\qquad \rightsquigarrow I_{\mathrm e}=\emptyset, \hat{\mathbf P}^{[12]}_{1\ldots5}=\begin{bmatrix}-2&0&2&0&-2\end{bmatrix}$

      and in this case $\Xi(A_1,A_2)=-A_2A_1$
\item $1[12]\typIV,4[12]\typIII \,\&\,2[12]\typIII,3[12]\typIV\,\&\,5[12]\typI$
      $\qquad \rightsquigarrow I_{\mathrm e}=\emptyset, \hat{\mathbf P}^{[12]}_{1\ldots5}=\begin{bmatrix}-2&0&-2&0&2\end{bmatrix}$

      and in this case $\Xi(A_1,A_2)=A_1^{-1}A_2^{-1}$
\item $1[12]\typII,4[12]\typI \,\&\,2[12]\typI,3[12]\typII\,\&\,5[12]\typI$
      $\qquad \rightsquigarrow I_{\mathrm e}=\emptyset, \hat{\mathbf P}^{[12]}_{1\ldots5}=\begin{bmatrix}0&2&0&2&2\end{bmatrix}$

      and in this case $\Xi(A_1,A_2)=-A_1A_2^{-1}$
\item $1[12]\typII,4[12]\typI \,\&\,2[12]\typIV,3[12]\typIII\,\&\,5[12]\typIV$
      $\qquad \rightsquigarrow I_{\mathrm e}=\emptyset, \hat{\mathbf P}^{[12]}_{1\ldots5}=\begin{bmatrix}0&-2&0&2&-2\end{bmatrix}$

      and in this case $\Xi(A_1,A_2)=A_2^{-1} A_1 $
\item $1[12]\typIII,4[12]\typIV \,\&\,2[12]\typI,3[12]\typII\,\&\,5[12]\typIV$
      $\qquad \rightsquigarrow I_{\mathrm e}=\emptyset, \hat{\mathbf P}^{[12]}_{1\ldots5}=\begin{bmatrix}0&2&0&-2&-2\end{bmatrix}$

      and in this case $\Xi(A_1,A_2)=A_2A_1^{-1}$
\item $1[12]\typIII,4[12]\typIV \,\&\,2[12]\typIV,3[12]\typIII\,\&\,5[12]\typI$
      $\qquad \rightsquigarrow I_{\mathrm e}=\emptyset, \hat{\mathbf P}^{[12]}_{1\ldots5}=\begin{bmatrix}0&-2&0&-2&2\end{bmatrix}$

      and in this case $\Xi(A_1,A_2)=-A_1^{-1}A_2$
\end{enumerate}

Remark: Conditions $5[12]\typII$,  $5[12]\typIII$,  $5[12]\typV$, $2[12]\langle+1\rangle$, $2[12]\langle-1\rangle$, $3[12]\langle-1\rangle$, $4[12]\langle+1\rangle$, $4[12]\langle-1\rangle$ are inconsistent; $3[12]\langle+1\rangle$ is trivial.

We can say that
conditions (1--12) are the primary conditions, (13--24) are composite conditions, (25--32) are extremal conditions.
\end{state}

\begin{state}\label{pan0}
For orthogonal invariant pseudoscalar FQ operations  $\Xi$, the consistent constellations of principal hyperscaling properties are as follows:
\begin{enumerate}
\item $1[12]\typV\Leftrightarrow 3[12]\typV$\label{pan1}
      $\qquad \rightsquigarrow I_{\mathrm{ne}}=\{2,4,5\}, \hat{\mathbf P}^{[12]}_{1\ldots5}=\begin{bmatrix}0&\hat q_2&0&0&0\end{bmatrix}$
\item $2[12]\typI \Leftrightarrow 3[12]\typII$\label{pan5}
      $\qquad \rightsquigarrow I_{\mathrm{ne}}=\{1,4,5\}, \hat{\mathbf P}^{[12]}_{1\ldots5}=\begin{bmatrix}0&2&0&0&0\end{bmatrix}$
\item $2[12]\typIV \Leftrightarrow 3[12]\typIII$\label{pan6}
      $\qquad \rightsquigarrow I_{\mathrm{ne}}=\{1,4,5\}, \hat{\mathbf P}^{[12]}_{1\ldots5}=\begin{bmatrix}0&-2&0&0&0\end{bmatrix}$
\item $2[12]\typII \Leftrightarrow 3[12]\typI$\label{pan3}
      $\qquad \rightsquigarrow I_{\mathrm{ne}}=\{1,4,5\}, \hat{\mathbf P}^{[12]}_{1\ldots5}=\begin{bmatrix}0&0&2&0&0\end{bmatrix}$
\item $2[12]\typIII \Leftrightarrow 3[12]\typIV$\label{pan4}
      $\qquad \rightsquigarrow I_{\mathrm{ne}}=\{1,4,5\}, \hat{\mathbf P}^{[12]}_{1\ldots5}=\begin{bmatrix}0&0&-2&0&0\end{bmatrix}$
\item $2[12]\typV\Leftrightarrow 4[12]\typV$\label{pan2}
      $\qquad \rightsquigarrow I_{\mathrm e}=\{3\}, \hat{\mathbf P}^{[12]}_{1\ldots5}=\begin{bmatrix}0&0&\hat q_3&0&0\end{bmatrix}$
\item $1[12]\typV,3[12]\typV \,\&\,2[12]\typI,3[12]\typII \Leftrightarrow6[12]\langle+1\rangle$
      $\quad \rightsquigarrow I_{\mathrm e}=\emptyset, \hat{\mathbf P}^{[12]}_{1\ldots5}=\begin{bmatrix}0&2&0&0&0\end{bmatrix}$

      and in this case $\Xi(A_1,A_2)=\mathcal A_L(A_1,A_2)$
\item $1[12]\typV,3[12]\typV \,\&\,2[12]\typIV,3[12]\typIII\Leftrightarrow6[12]\langle-1\rangle$
      $\quad \rightsquigarrow I_{\mathrm e}=\emptyset, \hat{\mathbf P}^{[12]}_{1\ldots5}=\begin{bmatrix}0&-2&0&0&0\end{bmatrix}$

      and in this case $\Xi(A_1,A_2)=\mathcal A_R(A_1,A_2)$
\item $2[12]\typII,3[12]\typI \,\&\,2[12]\typV,4[12]\typV\Leftrightarrow7[12]\langle+1\rangle$
      $\quad \rightsquigarrow I_{\mathrm e}=\emptyset, \hat{\mathbf P}^{[12]}_{1\ldots5}=\begin{bmatrix}0&0&2&0&0\end{bmatrix}$

      and in this case $\Xi(A_1,A_2)=\mathcal D(A_1,A_2)$

\item $2[12]\typIII,3[12]\typIV \,\&\,2[12]\typV,4[12]\typV\Leftrightarrow7[12]\langle-1\rangle$
      $\quad \rightsquigarrow I_{\mathrm e}=\emptyset, \hat{\mathbf P}^{[12]}_{1\ldots5}=\begin{bmatrix}0&0&-2&0&0\end{bmatrix}$

      and in this case $\Xi(A_1,A_2)=-\mathcal D(A_1,A_2)^{-1}$
\end{enumerate}

We can say that conditions (1--6) are the primary conditions, (7--10) are extremal conditions.
\end{state}
\begin{state}For orthogonal invariant pseudoscalar FQ operations,
\begin{itemize}
\item an example for \ref{pan0}\eqref{pan1} is $\Xi=\mathcal A_C$ with
$\hat q_2=0$, $\hat{\mathbf P}^{[12]}_{1\ldots5}=\begin{bmatrix}0&0&0&0&0\end{bmatrix}$;
\item   examples for \ref{pan0}\eqref{pan5} are $\Xi=\mathcal D\mathcal L^{-1},\pol \mathcal D\mathcal L^{-1},-\mathcal L\mathcal D^{-1},$
$\mathcal D^{\mathrm{x}}(\mathcal L^{\mathrm{x}})^{-1}=-\mathcal L^{\mathrm{cx}}(\mathcal D^{\mathrm{cx}})^{-1},$
$\mathcal D^{\mathrm{cx}}(\mathcal L^{\mathrm{cx}})^{-1}=-\mathcal L^{\mathrm{x}}(\mathcal D^{\mathrm{x}})^{-1}$;
\item   examples for \ref{pan0}\eqref{pan6} are $\Xi=\mathcal L^{-1}\mathcal D,\pol\mathcal L^{-1}\mathcal D,-\mathcal D^{-1}\mathcal L,$
$(\mathcal L^{\mathrm{x}})^{-1}\mathcal D^{\mathrm{x}}=-(\mathcal D^{\mathrm{cx}})^{-1}\mathcal L^{\mathrm{cx}},$
$(\mathcal L^{\mathrm{cx}})^{-1}\mathcal D^{\mathrm{cx}}=-(\mathcal D^{\mathrm{x}})^{-1}\mathcal L^{\mathrm{x}}$;
\item   examples for \ref{pan0}\eqref{pan3} are $\Xi=\mathcal D^{\mathrm{fx}},\mathcal D^{\mathrm{x}},\mathcal D^{\mathrm{cx}}$;
\item   examples for \ref{pan0}\eqref{pan4} are $\Xi=-(\mathcal D^{\mathrm{fx}})^{-1},-(\mathcal D^{\mathrm{x}})^{-1},-(\mathcal D^{\mathrm{cx}})^{-1}$;
\item an  example for \ref{pan0}\eqref{pan2} is $\Xi=\mathcal A_D=\pol \mathcal D$ with $\hat q_3=0$,
$\hat{\mathbf P}^{[12]}_{1\ldots5}=\begin{bmatrix}0&0&0&0&0\end{bmatrix}$.
\end{itemize}
\end{state}
What we see is that only  very few FQ operations are characterized by Clifford conservativity,
principal hyperscalings and orthogonal invariance alone.
The situation improves if one allows to combine them with scalar scaling conditions.
For example, in case of \ref{san0}\eqref{san2}, after only  the index $3$ left, a simple
scalar homogeneity property (i. e. scalar scaling in $\hat r_3$) with $\hat p_3=\alpha$ is sufficient to fix the FQ operation.

The statements above, in this form, are, of course, conjectural, and their proofs should be somewhat
longish due, if not else, to the length of the statements themselves.
However, certain restrictive aspects of them (like some inconsistencies) can be checked rather easily.
The general picture they suggest is that hyperscaling conditions are heterogeneous, but they do not
describe operations with any very specific properties but the arithmetically very simplest ones.
In particular, hyperscaling conditions already limit first order behaviour severely
(especially in the presence of orthogonal and transposition invariance properties),
which makes them unsuitable for certain classes of operations.
\begin{point}
Alternatively, one can try the combined hyperscaling conditions in $\tilde r_4$ and $\tilde r_5$.
In this case $L=0$ can be assumed anyway. $\tilde 4[s](J,0,\alpha,\beta)$ and
$\tilde 5[s](J,0,\alpha,\beta)$
reads as (with $\pm=+$ and $\pm=-$, respectively)
\begin{align}
\tfrac12\hat p^{[s]}_{4}\pm\tfrac12\hat p^{[s]}_{5}&=(\alpha+\beta)\hat p^{[s]}\notag\\
\tfrac12\hat p^{[s]}_{4,j,\cdots}\pm\tfrac12\hat p^{[s]}_{5,j,\cdots} &=
\alpha \hat p^{[s]}_{j,\cdots}
-\tfrac14(J\pm1)\hat p^{[s]}_{5*j,\cdots}\mp\tfrac14(J\pm1)\hat p^{[s]}_{4*j,\cdots}\notag\\
\tfrac12\hat p^{[s]}_{\cdots,i,4}\pm\tfrac12\hat p^{[s]}_{\cdots,i,5}&=
\tfrac14(J\mp1)\hat p^{[s]}_{\cdots,i*5}\pm\tfrac14(J\mp1)\hat p^{[s]}_{\cdots,i*4}
+\beta\hat p^{[s]}_{\cdots,i}\notag\\
\tfrac12\hat p^{[s]}_{\cdots,i,4,j,\cdots}\pm\tfrac12\hat p^{[s]}_{\cdots,i,5,j,\cdots}&=
\tfrac14(J\mp1)\hat p^{[s]}_{\cdots,i*5,j,\cdots}\pm\tfrac14(J\mp1)\hat p^{[s]}_{\cdots,i*4,j,\cdots}
\notag\\
&-\tfrac14(J\pm1)\hat p^{[s]}_{\cdots,i,5*j,\cdots}\mp\tfrac14(J\pm1)\hat p^{[s]}_{\cdots,i,4*j,\cdots}\qquad.
\notag
\end{align}
Similarly, as before, they lead to various principal hyperscaling conditions, but they offer little more in the orthogonal invariant case:
\begin{state}If an scalar or symmetric vectorial or pseudoscalar orthogonal invariant  FQ operation $\Xi$
satisfies a  hyperscaling condition in $\tilde r_4$ or $\tilde r_5$, then it is one of the following:
\begin{enumerate}
\item $\Xi=1$ (satisfies, for example,  $\tilde 4,\tilde 5[0]\typII,\typIII,\typV$)
\item $\Xi=\Id$ (satisfies, for example,  $\tilde 4[1]\typVI,\typVIII$)
\item $\Xi=\mathcal T_ {RR}$ (satisfies  $\tilde 4[1]\typVII$)
\item $\Xi=\mathcal T_ {LL}$ (satisfies  $\tilde 4[1]\typIX$)
\item $\Xi=\mathcal O^{\mathrm{afSy}}$ (satisfies, for example,  $\tilde 4[1]\typV$)
\item $\Xi=\hat q_5\Id+(1-\hat q_5)\mathcal O^{\mathrm{fSy}}$ (satisfies, for example,  $\tilde 5[1]\typV$)
\item $\Xi=A_R$ (satisfies, for example,  $\tilde 4[12]\typV$)
\item $\Xi=A_L$ (satisfies, for example,  $\tilde 5[12]\typV$)
\end{enumerate}
\end{state}
One can also try hyperscaling condition with $L\neq0$.
This leads to some principal types with $L=\pm1$ and a more complicated situation, but not much new in regard of orthogonal invariant FQ operations.
\end{point}
\begin{point}
Or, we can carry out the same computations  in the monoaxial regime.
Then we deal only with variables $\hat r_3,\hat r_4,\hat r_5$ but the same principal types can be used.
In fact, what happens is that we get the hyperscaling conditions in a much cleaner form as the principal
hyperscaling conditions do not ``glue'' together as before.
Nevertheless, ``interactions'' between them are possible if  more of them are imposed.
But even after the axial extensions we do not really arrive to essentially new examples compared to what we have seen.
\end{point}
So, while scalar scalings are much weaker, they can be used more flexibly  than hyperscalings.
Also, processes like axial extensions produce similar reductions but less restrictive.

\end{document}